\documentclass [11pt,oneside]{article}

\reversemarginpar \textwidth=14cm \textheight=19cm

\usepackage{amssymb}
\usepackage{amsfonts}
\usepackage{amsmath}
\usepackage{amsthm}
\usepackage[dvips]{graphicx}

\newtheorem{lemma}{Lemma}[section] \newtheorem{teo}[lemma]{Theorem}
\newtheorem{rem}[lemma]{Remark} \newtheorem{prop}[lemma]{Proposition}
\newtheorem{cor}[lemma]{Corollary}

\newcommand{\matN}{\ensuremath {\mathbb{N}}}
\newcommand{\matR} {\ensuremath {\mathbb{R}}}

\newcommand{\matZ} {\ensuremath {\mathbb{Z}}}

\newcommand{\matRP} {\ensuremath {\mathbb{RP}}}

\newcommand{\calX} {\ensuremath {\mathcal{X}}}
\newcommand{\calR} {\ensuremath {\mathcal{R}}}
\newcommand{\calN} {\ensuremath {\mathcal{N}}}
\newcommand{\calK} {\ensuremath {\mathcal{K}}}

\newcommand{\calL} {\ensuremath {\mathcal{L}}}

\newcommand{\calB} {\ensuremath {\mathcal{B}}}

\newcommand{\calT} {\ensuremath {\mathcal{T}}}

\newcommand{\calH}{\ensuremath {\mathcal{H}}}

\newcommand{\nota} [1] {\caption{\footnotesize{#1}}}

\def\interior#1{{\rm int}(#1)}

\font\titsc=cmcsc10 scaled 1200

\author{Bruno \titsc{Martelli} \and Carlo \titsc{Petronio}}

\title{3-manifolds having complexity at most $9$}

\begin{document}

\maketitle

\tableofcontents

\section{Introduction} \label{introduction}

This paper is devoted to the theoretical description and illustration
of results of an algorithm which has enabled us to give a complete
list, without repetitions, of all closed oriented irreducible
$3$-manifolds of complexity up to $9$. More interestingly, we have
actually been able to give a ``name'' to each such manifold, {\em
i.e.}~to recognize its canonical decomposition into Seifert fibered
spaces and hyperbolic manifolds already considered by other authors.
The complexity we are referring to here is that introduced by Matveev
(\cite{Mat}, see also \cite{MaFo}), given by the minimal number of
vertices of a simple spine (this has been proved in \cite{Mat} to be
equal to the minimal number of tetrahedra in a triangulation).

Our algorithm relies on a structural result on closed 3-manifolds.
Namely, we show that all closed 3-manifolds can be obtained by
combining, in a suitable sense, building blocks taken from a certain
list which, at least up to complexity 9, is dramatically shorter than
the list of all manifolds. The building blocks are called
\emph{bricks}, they are bounded by tori, and these tori carry a
``marking'' given by an embedded trivalent graph. Moreover, the
combination of two bricks corresponds to the identification of two
boundary tori. The main definitions and results of the theory of
decomposition into bricks are stated in the rest of the present
introduction and proved in the body of the paper.

Before turning to bricks, let us mention the most interesting experimental
results about complexity $9$ which our algorithm has allowed us to
discover. Recall first that it was already known to Matveev \cite{Mat}
that up to complexity $8$ all manifolds are graph-manifolds; tables up to
complexity $6$ are in \cite{Mat2}, and up to $7$ in \cite{Ov}. Now, we can
show that there are $1156$ manifolds of complexity $9$, $272$ of them are
lens spaces, $863$ are more general graph-manifolds which do not contain
non-separating tori, $17$ of them are torus bundles over $S^1$, $10$ of
them are graph-manifolds with graph \begin{picture}(21,8)
\put(4.2,4){\circle{8}} \put(8,4){\circle*{2}} \put(8,4){\line(1,0){8}}
\put(16,4){\circle*{2}} \end{picture}, and there are also $4$ closed
hyperbolic manifolds. More importantly, these $4$ manifolds turn out to be
precisely those of least known volume~\cite{Ho-We}, in accordance
with the ideas about complexity and volume stated in \cite{MaFo}.

\subsection{Bricks and assemblings of bricks} Throughout this paper we
will work in the PL category, and by \emph{manifold} we will always mean a
compact orientable 3-manifold, possibly with boundary. We will call
\emph{triod} the graph with two vertices and three edges all joining one
vertex to the other one. Note that a triod $\theta$ can be embedded in a
torus $T$ so that $T\setminus\theta$ is an open 2-disc. A pair $(M,X)$ is
said to be a \emph{manifold with triods} if $M$ is a manifold with
boundary consisting of tori $T_1, \ldots, T_n$ and $X$ is a set of triods
$\{\theta_1, \ldots, \theta_n\}$, with $\theta_i$ embedded in $T_i$ so
that $T_i \setminus \theta_i$ is a disc. The case where $n=0$ and
$X=\emptyset$, so $M$ is closed, is admitted.

Let $\calX$ be the set of all manifolds with triods (up to equivalence
induced by homeomorphism of manifolds). If $M$ has non-empty boundary
consisting of tori, then there are infinitely many inequivalent ways to
embed triods in these tori, so there are infinitely many inequivalent
pairs $(M,X)$ based on the same $M$. On the contrary, if $M$ is closed,
then there is a unique element $(M,\emptyset)\in\calX$ based on $M$.
Therefore the set of all closed orientable manifolds can be viewed as a
subset of $\calX$.

We will now describe three operations on $\calX$ and state the crucial
properties of a complexity function on $\calX$ introduced and discussed in
detail below in Section~\ref{ct}.

\paragraph{Connected sum.} The operation of connected sum ``far from the
boundary'' obviously extends from manifolds to manifolds with triods.
Namely, given $(M,X)$ and $(M',X')$ in $\calX$, we define $(M,X)\#(M',X')$
as $(M\# M',X\cup X')$, where $M\# M'$ is one of the two possible
connected sums of $M$ and $M'$ (recall that our manifolds are orientable
but not oriented). Of course $(S^3,\emptyset)\in\calX$ is the identity
element for operation $\#$. We will call a pair $(M,X)$ \emph{prime} if
$M$ is, \emph{i.e.}~if $(M,X)$ cannot be expressed as a connected sum of
pairs different from $(S^3,\emptyset)$.

\paragraph{Assembling.} Given $(M,X)$ and $(M',X')$ in $\calX$, we pick
triods $\theta_i \in X$ and $\theta'_{i'} \in X'$ and choose a
homeomorphism $\psi:T_i \to T_{i'}'$ such that $\psi(\theta_i) =
\theta'_{i'}$. We can then construct the manifold with triods $(N,Y) = (M
\cup_{\psi} M',(X\cup X')\setminus \{\theta_i, \theta_{i'}'\})$. We call
this an \emph{assembling} of $(M,X)$ and $(M',X')$ and we write
$(N,Y)=(M,X)\oplus(M',X')$. Of course two given elements of $\calX$ can
only be assembled in a finite number of inequivalent ways.

Operation $\oplus$ has an identity element, and in a special case
it is the inverse operation of $\#$. Below we will need to exclude
these types of assembling, so we describe them in detail.
First, set $B_0=(T\times [0,1],\{\theta \times \{0\}, \theta \times \{1\}\})$,
where $T$ is the torus and $\theta \subset T$ is a triod such that
$T\setminus \theta$ is a disc ($B_0$ is well-defined up to equivalence).
Of course if we assemble any $(M,X)\in\calX$ with $B_0$ we get $(M,X)$ again.

Let $H$ be the solid torus and let $(H,\{\theta\})$ and
$(H,\{\theta'\})$ be elements of $\calX$ based on $H$. Assume that
there exists a homeomorphism $\partial H \to
\partial H$ with $\psi(\theta)=\theta'$  such that
$(H,\{\theta\})\oplus(H,\{\theta'\})$ performed along $\psi$ gives
$(S^3,\emptyset)$ as a result. Then for any $(M,X)\in\calX$ we have
$((M,X)\#(H,\{\theta\}))\oplus(H,\{\theta'\})=(M,X)$ if we use the same
$\psi$.

This discussion motivates the following definition.
An assembling $(M,X)\oplus(M',X')$ is called \emph{trivial} if, up to
interchanging $(M,X)$ and $(M',X')$, one of the following holds:
\begin{itemize}
\item $(M',X')=B_0$, or
\item $(M',X')=(H,\{\theta'\})$ is a solid torus with triod,
and $(M,X)$ can be decomposed as $(M,X)=(N,Y)\# (H,\{\theta\})$
so that $(N,Y)\ne (S^3,\emptyset)$ and
the assembling identifies $\theta$ to $\theta'$ and
$(H,\{\theta\})\oplus(H,\{\theta'\})=(S^3,\emptyset)$.
\end{itemize}

\paragraph{Self-assembling.} Given $(M,X)\in\calX$, we pick two distinct
triods $\theta_i, \theta_{i'} \in X$, we choose a homeomorphism $\psi:T_i
\to T_{i'}$ such that $\psi(\theta_i)$ and $\theta_{i'}$ intersect
transversely in two points, and we construct the manifold with triods
$(N,Y) = (M_{\psi}, X \setminus \{\theta_i, \theta_{i'}\})$. We call this
a \emph{self-assembling} of $(M,X)$ and we write $(N,Y) = \odot(M,X)$. As
above, only a finite number of self-assemblings of a given element of
$\calX$ are possible.

In the sequel it will be convenient to refer to a combination of
assemblings and self-assemblings of pairs just as an \emph{assembling}.
Note that of course we can do the assemblings first and the
self-assemblings in the end.

\paragraph{A complexity on $\calX$.} One of the main ingredients of the
present paper is the extension of Matveev's definition of
complexity~\cite{Mat} from closed manifolds to manifolds with triods.
We warn the reader that Matveev's complexity $c(M)$ is defined also
when $\partial M\neq\emptyset$, but our definition will be different
in this case, namely we will have $c(M,X)=c(M)$ only when
$X=\emptyset$, \emph{i.e.}~when $M$ is closed. The key properties of
$c$, proved below, are additivity with respect to connected sum and
subadditivity with respect to assembling. More precisely, we will
construct in Subsection~\ref{ct-def} a function $c:\calX \to \matN$
and show in Subsection~\ref{ct-properties} that it enjoys the
following properties:

\begin{enumerate}

\item\label{extension:property} $c(M,\emptyset)=c(M)$ for any
$(M,\emptyset)\in\calX$;

\item\label{con:sum:add:property} $c((M,X)\#(M',X'))=c(M,X)+c(M',X')$;

\item\label{assemble:subad:property} $c((M,X)\oplus(M',X'))\le
c(M,X)+c(M',X')$. Moreover, when equality holds and the assembling is non-trivial,
we have that $(M,X)\oplus(M',X')$ is prime if and only if both $(M,X)$ and $(M',X')$
are;

\item\label{self-assemble:subad:property} $c(\odot(M,X))\le c(M,X)+6$.
Moreover, when equality holds, we have that $\odot(M,X)$ is prime if and only
if $(M,X)$ is;

\item\label{finite:property} for any $n\geq 0$ there is only a finite
number of prime pairs $(M,X)\in\calX$ with $c(M,X)\le n$.

\end{enumerate}

Now let $\calX^{\rm pr}\subset\calX$ be the set consisting of prime
pairs. An assembling is called \emph{sharp} if it is non-trivial and
the inequality of point~(\ref{assemble:subad:property}) above is
actually an equality. Similarly, a self-assembling is \emph{sharp} if
in~(\ref{self-assemble:subad:property}) we have an equality. We will
say that a prime pair $(M,X)\in\calX^{\rm pr}$ is
a \emph{brick} if it cannot be expressed as the result of a sharp
assembling or a sharp self-assembling. The following easy result will
be proved in Subsection~\ref{ct-def} (one could actually also deduce
it from property~(\ref{finite:property}), but we will refrain from
doing this):

\begin{lemma}\label{zero:compl:lem} The pair $B_0$ is the only
$(M,X)\in\calX$ such that $c(M,X)=0$ and $X$ contains at least two triods.
\end{lemma}

Induction on complexity now readily implies the following:

\begin{prop} \label{main-prop} Every prime manifold with triods can be
obtained as a sharp-assembling of some bricks. \end{prop}

We define now $\calB\subset\calX^{\rm pr}$ as the set of all bricks, and
note that $\calB$ naturally splits as $\calB^0 \sqcup \calB^1$, where
$\calB^0$ is the set of all $(M,X)\in \calB$ with $X=\emptyset$
(\emph{i.e.}~$M$ is closed). Pairs in $\calB^0$ cannot be used for an
assembling or self-assembling, since they have no boundary. Let $\calB_n^j
\subset \calB^j$, for $j=0,1$, and $\calX_n \subset \calX$ be the subsets
consisting of pairs having complexity $n$. Proposition~\ref{main-prop} and
the properties of $c$ stated above now imply that $$\calX^{\rm pr}_{\le n}
=\calB^0_{\le n}\; \cup\; \left\{\odot^k(B^1\oplus \ldots \oplus B^h):\
B^i\in\calB^1_{\le n},\ \sum c(B^i)+6k \le n\right\}.$$

If one can give an unambiguous name to each closed $\odot^k(B^1\oplus
\ldots \oplus B^h)$, then the set of all closed prime manifolds having
complexity at most $n$ is easily constructed from $\calB_{\le n}$ by
listing the (finite number of) closed manifolds obtained in this way, and
by then removing duplicates. For $n\le 9$ it turns out that $\calB_{\le
n}$ consists of a very few atoroidal manifolds (with triods), and it is
experimentally not so hard to give a name to each closed manifold of the
form $\odot^k(B^1\oplus \ldots \oplus B^h)$. We will provide more details
below on the recognition issue (after listing the bricks explicitly), but
we want to emphasize here that the vast majority of computer time in the
implementation of our algorithm was taken by the determination of bricks.
Taking the list of bricks for granted, the reader could with some patience
reproduce the list of manifolds by himself.

\subsection{Bricks and manifolds up to complexity 9}\label{results} The
algorithm which will be explained in Section~\ref{algorithm} has enabled
us to explicitly find $\calB_{\le 9}^0$ and $\calB_{\le9}^1$. The former
consists of 19 closed manifolds naturally coming in two families $C_{i,j}$
and $E_k$, and the latter consists of only 11 manifolds with triods,
denoted by $B_0,\dots,B_{10}$ (where $B_0$ is the same as defined above).
The elements of $\calB_{\le9}^0$ are all Seifert fibered over $S^2$ with 3
exceptional fibers. In order to describe the elements of $\calB_{\le9}^1$
we need a way to encode the possible ways a triod can sit in a torus.

\begin{rem} \emph{Let $T$ be a torus. Let $\calT$ be the set of unordered
triples $\{a,b,c\}$ of elements of $H_1(T)$, such that every pair of
elements in $\{a,b,c\}$ generates $H_1(T)$, and $a+b+c=0$. Let $\theta
\subset T$ be a triod such that $T\setminus \theta$ is a disc: inside
$\theta$ we can find $3$ distinct closed curves, which can be oriented in
order to form a triple $\{a,b,c\}\in\calT$. The only two triples we can
get like this are $\{a,b,c\}$ and $\{-a,-b,-c\}$. Conversely, each triple
$\{a,b,c\}\in\calT$ determines a triod $\theta \subset T$. It follows that
triods (up to isotopy) are in one-to-one correspondence with elements of
$\calT/\matZ_2$, where the non-trivial element of $\matZ_2$ acts mapping
$\{a,b,c\}$ to $\{-a,-b,-c\}$.} \end{rem}

\paragraph{Bricks.} In Table~\ref{brick:table} we list the elements in
$\calB_{\le 9}^1$, as produced by our algorithm, where $D_k$ is the
disc with $k$ holes, and the usual notation for Seifert manifolds and
cusped hyperbolic manifolds~\cite{Ca-Hi-We} is employed. Note that
$c(M)\neq c(M,X)$ is the complexity of $M$ in the usual sense
\cite{Mat}, defined for any compact 3-manifold. Every $M$ turns out
to be atoroidal. In order to describe triods as elements in $\calT$,
we must fix a basis $(\mu_i, \lambda_i)$ for $H_1(T_i)$ for each
$T_i$ in $\partial M$. When $M$ is Seifert, by removing fibered
neighbourhoods of the exceptional fibers we get $D_k \times S^1$ with
the product fibration. Then we choose $\lambda_i$ to be a fiber and
$\mu_i$ to be a component of $\partial D_k\times\{{\rm point}\}$,
with orientations chosen so that $(\mu_i, \lambda_i)$ is a positively
oriented basis. When $M$ is hyperbolic we choose $\mu_i$ and
$\lambda_i$ to be respectively the first and second shortest
geodesic, with orientations such that $(\mu_i, \lambda_i)$ is a
positively oriented basis. In both cases, taking $(-\mu_i,
-\lambda_i)$ instead of $(\mu_i,\lambda_i)$ as a basis does not make
any difference, since triples in $\calT$ are defined up to sign.

\begin{table}
\begin{center}\begin{tabular} {|c|c|c|c|c|}
\hline $(M,X)$ & $c(M,X)$ & $M$ & $X$ & $c(M)$ \\ \hline \hline $B_0$
& $0$ & $D_1 \times S^1$ & $\{(1,0),(0,1),(-1,-1))\}$ & $0$ \\ & & &
$\{(1,0),(0,-1),(-1,1)\}$ & \\ \hline $B_1$ & $0$ & $D_0 \times S^1$
& $\{(1,0),(0,1),(-1,-1)\}$ & $0$
\\ \hline
$B_2$ & $0$ & $D_0 \times S^1$ & $\{(0,1),(1,1),(-1,-2)\}$ & $0$
\\ \hline
$B_3$ & $1$ & $D_1 \times S^1$ & $\{(1,0),(0,1),(-1,-1)\}$ & $0$
\\ & & & $\{(1,0),(0,1),(-1,-1)\}$ & \\ \hline
$B_4$ & $3$ & $D_2 \times S^1$ & $\{(1,0),(0,1),(-1,-1)\}$ & $0$
\\ & & & $\{(1,0),(0,1),(-1,-1)\}$ &
\\ & & & $\{(1,0),(0,-1),(-1,1)\}$ & \\ \hline
$B_5$ & $8$ & $(D_0,(2,1),(3,1))$ & $\{(1,-1),(5,-4),(-6,5)\}$ & $0$
\\ \hline $B_6$ & $8$ & $M2_2^1$ & $\{(1,0),(0,-1),(-1,1)\}$ &
$2$ \\ \hline $B_7$ & $9$ & $M3_4^1$ & $\{(1,0),(0,-1),(-1,1)\}$ &
$3$ \\ \hline $B_8$ & $9$ & $M4^2_1$ & $\{(1,0),(0,1),(-1,-1)\}$ &
$4$ \\ & & & $\{(1,0),(0,1),(-1,-1)\}$ &
\\ \hline $B_9$ & $9$ & $M6^3_1$ & $\{(1,0),(0,-1),(-1,1)\}$ &
$6$ \\ & & & $\{(1,0),(0,-1),(-1,1)\}$ &
\\ & & & $\{(1,0),(0,-1),(-1,1)\}$ & \\ \hline $B_{10}$ & $9$ & $M6^3_1$ &
$\{(1,0),(0,-1),(-1,1)\}$ & $6$ \\ & & & $\{(1,0),(0,-1),(-1,1)\}$ &
\\ & & & $\{(1,0),(0,1),(-1,-1)\}$ & \\ \hline
\end{tabular}
\nota{Bricks up to complexity 9.}  \label{brick:table}
\end{center}
\end{table}

\paragraph{From bricks to manifolds.} As pointed out in the introduction,
a list of all closed orientable prime manifolds with complexity at most
$9$ can be compiled by listing and recognizing all closed manifolds
obtained by assembling bricks $B^1,\dots,B^h$ of $\calB_{\le9}$ and then
self-assembling $k$ times, with $\sum c(B^i)+6k\le 9$. We explain here the
points which make this listing and recognition feasible. Note first that,
by the bound on complexity, only a few assemblings, and no
self-assembling, will involve $B_5,\ldots,B_{10}$. We also know that $B_0$
must not be used for assemblings. Moreover we can eliminate from the list
all assemblings which we know \emph{a priori} not to be sharp. For
instance we have the following (proved in Section~\ref{ct}):

\begin{prop} \label{not:sharp} If $(M,X)\in\calX^{\rm pr}$
and $(M,X)\oplus B_1$ is sharp, then $(M,X)$ is either $B_1$
or $B_2$.\end{prop}

\begin{figure} \begin{center}
\includegraphics{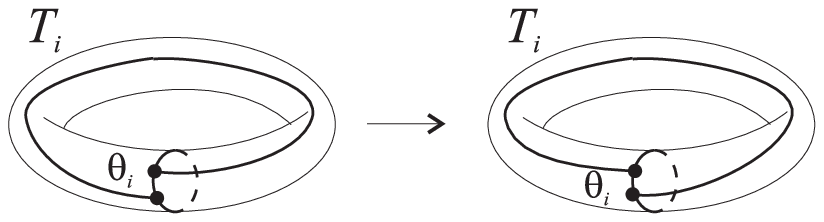}
\nota{The effect of an assembling with $B_3$.} \label{b3eff}
\end{center} \end{figure}

Concerning recognition, we note now that the effect of assembling $B_2$ or
$B_3$ is very easy to describe. Since $B_2$ is a solid torus, the
assembling with $B_2$ along some boundary component $T_i$ corresponds to a
Dehn filling of $T_i$. Finitely many different fillings are possible, and
they are determined by the position in $T_i$ of the triod $\theta_i$. Now
$B_3\cong (T\times[0,1], \{\theta_0\times\{0\},\theta_1\times\{1\}\})$
with $\theta_0\neq\theta_1$. (Even if in Table~\ref{brick:table} the
triples describing the triods are the same, the triods are not the same,
because they lie on different boundary components, so the bases of
homology are different due to orientation.) More precisely, the assembling
with $B_3$ along $T_i$ corresponds to changing the position of the triod
$\theta_i$ as in Fig.~\ref{b3eff}. Summing up, the successive assembling
along $T_i$ of some $B_3$'s followed by the assembling of one $B_2$ still
corresponds to a Dehn filling of $T_i$. One actually sees that all Dehn
fillings can be generated like this, but of course the bound on complexity
allows to consider only finitely many of them.

Turning to $B_4$ and $B_5$, we note that they naturally come with a
Seifert fibered structure, so any manifold generated by $B_2,\ldots,B_5$
is a graph manifold, whose graph and gluing matrices are readily deduced
from the pattern of assemblings giving the manifold. Since there are
algorithms checking whether two such set of data give the same manifold,
recognition is not a problem at this level.

Getting to assemblings involving $B_6,\ldots,B_{10}$, one first notes that
they can only be assembled with $B_2$ and $B_3$, and not in many ways.
Next, one checks by direct comparison with the tables
in~\cite{Ho-We} that 4 of the resulting manifolds are the 4
hyperbolic closed manifolds with least known volume. The following fact
(proven in Section~\ref{concluding}) concludes our investigation:

\begin{prop} \label{not:sharp2} Let
$M$ be a closed manifold with $c(M)\le 9$ obtained by
assembling a brick in $\{B_6\ldots,B_{10}\}$ and some $B_2$'s and $B_3$'s.
Then either $M$ is one of the 4 hyperbolic manifolds just described, or
the assembling is not sharp. \end{prop}

\paragraph{Manifolds.} Table~\ref{manif:table} contains the data our
algorithm has allowed us to discover about closed orientable prime
manifolds having complexity $c$ for $c\le 9$. We have divided the
manifolds into three groups, given respectively by the elements which may
be obtained by sharp-assembling $B_0,\dots,B_{10}$ but without
self-assembling, by those which require a self-assembling, and by those of
$\calB^0$. The three groups have been further split to give a more precise
idea of which bricks are needed to generate a manifold: in particular, the
vast majority of manifolds (with 11 exceptions out of 1156 manifolds in
complexity 9) are obtained assembling $\{B_2,B_3,B_4\}$, and only a few
manifolds actually require a self-assembling. An important convention in
the table is that manifolds already considered in a certain line are not
considered again in subsequent lines: some manifolds can be split into
bricks in distinct ways.

\begin{table}
\begin{center} \begin{tabular}{|l|c|c|c|c|c|c|c|c|c|c|}
\hline Vertices & $0$ & $1$ & $2$ & $3$ & $4$ & $5$ & $6$ & $7$ & $8$ & $9$
\\ \hline\hline $\langle B_1 \rangle_{\rm non-self}$ & $2$ & & & & & & & & & \\
$\langle B_2 \rangle_{\rm non-self}$ & $2$ & & & & & & & & & \\
$\langle B_2,B_3 \rangle_{\rm non-self}$ & & $2$ & $3$ & $6$ & $10$ &
$20$ & $36$ & $72$ & $136$ & $272$ \\ $\langle B_2,B_3, B_4
\rangle_{\rm non-self}$ & & & & & $2$ & $8$ & $32$ & $97$ & $292$ &
$856$ \\ $\langle B_2, B_3, B_5 \rangle_{\rm non-self}$ & & & & & & &
& & $1$ & $3$ \\ $\langle B_2,B_3, B_6 \rangle_{\rm non-self}$ & & &
& & & & & & & $2$ \\ $\langle B_2,B_7\rangle_{\rm non-self}$ & & & &
& & & & & & $1$ \\ $\langle B_2, B_8\rangle_{\rm non-self}$ & & & & &
& & & & & $1$
\\ \hline $\langle B_0\rangle_{\rm self}$ & & & & & & & $5$ & &
& \\ $\langle B_3\rangle_{\rm self}$ & & & & & & & & $3$ & $3$ & $7$
\\ $\langle B_2, B_4\rangle_{\rm self}$ & & & & & & & & & & $10$
\\ \hline $C_{i,j}$ & & & $1$ & $1$ & $2$ & $2$ & $1$ & $2$ & $3$ & $3$
\\ $E_k$ & & & & & & $1$ & $0$ & $1$ & $1$ & $1$ \\ \hline \hline
Total & $4$ & $2$ & $4$ & $7$ & $14$ & $31$ & $74$ & $175$ & $436$ &
$1156$
\\ \hline
\end{tabular}
\nota{Manifolds up to complexity 9.} \label{manif:table}
\end{center}
\end{table}

It follows from the topology of the bricks that the elements of $\langle
B_2,B_3\rangle_{\rm non-self}$ are all lens spaces, those of $\langle
B_2,B_3,B_4 \rangle_{\rm non-self}$ and $\langle B_2,B_3,B_5 \rangle_{\rm
non-self}$ are more general graph-manifolds whose graph is a tree, those
of $\langle B_3 \rangle_{\rm self}$ are torus bundles over $S^1$ and those
of $\langle B_2, B_4\rangle_{\rm self}$ are graph-manifolds with graph
\begin{picture}(21,8) \put(4.2,4){\circle{8}} \put(8,4){\circle*{2}}
\put(8,4){\line(1,0){8}} \put(16,4){\circle*{2}} \end{picture}. As already
mentioned, and explained in detail below in Section~\ref{concluding}, the
elements of $\calB^0$ (namely the $C_{i,j}$'s and $E_k$'s) are all Seifert
fibered over $S^2$ with 3 exceptional fibers.

\section{The complexity function} \label{ct} In this section we extend
Matveev's complexity~\cite{Mat} to manifolds with triods, and we state and
prove its properties.

\subsection{Definition of complexity} \label{ct-def}

A compact polyhedron $P$ is called \emph{simple} if the link of every
point of $P$ can be embedded in the space given by a circle with three
radii. The points having the whole of this space as a link are called
\emph{vertices}: they are isolated and therefore finite in number.

Let $(M,X)$ be a manifold with triods. A sub-polyhedron $P$ of $M$ is said
to be a \emph{skeleton} of the pair $(M,X)$ if

\begin{itemize}

\item $P\cup \partial M$ is simple, and $M\setminus (P\cup\partial M)$ is
an open ball;

\item $P\cap\partial M = X$.
\end{itemize}
Note that each open disc $T_i \setminus \theta_i$ is automatically
adjacent to the ball $M\setminus (P\cup\partial M)$, $P$ is simple, and
the vertices of $P$ cannot lie on $\partial M$. Note also that when $\#X
= 1$ then $P$ is a \emph{spine} of $M$ (\emph{i.e.}~$M$ collapses onto
$P$), and when $\#X=0$ (\emph{i.e.}~when $M$ is closed) then $P$ is a
spine in the usual sense~\cite{Mat}, namely $M\setminus\{{\rm point}\}$
collapses onto $P$. When $\#X \ge 2$, then $M$ does not collapse onto
$P$.

\begin{rem} \emph{It is easy to prove that every $(M,X)\in\calX$ has a
skeleton: take any simple spine $Q$ of $M\setminus\{{\rm point}\}$, so that
$M\setminus Q=\partial M \times [0,1)\cup B^3$, and assume that the
various $\theta_i\times[0,1)$'s are incident in a generic way to $Q$
and to each other (here of course the $\theta_i$'s are the triods in
$X$). Taking the union of $Q$ with the $\theta_i \times [0,1)$'s we
get a simple $Q'$ such that $M\setminus(Q'\cup\partial M)$ consists
of $\#X+1$ balls. Then we get a skeleton of $(M,X)$ by puncturing
$\#X$ suitably chosen 2-discs embedded in $Q'$, so to get one ball
only in the complement.}
\end{rem}

\begin{rem}{\em A definition of skeleton analogous to our one was given
in~\cite{Tu-Vi} for any compact manifold with any trivalent graph in its
boundary. The notion of complexity we will now introduce extends to any
such object.} \end{rem}

We say that a skeleton of $(M,X)$ is \emph{nuclear} if it does not
collapse to a subpolyhedron which is also a skeleton of $(M,X)$. We say
that a skeleton $P$ of $(M,X)\in\calX$ is \emph{minimal} if it is nuclear
and no other skeleton of $(M,X)$ has fewer vertices. We define now the
\emph{complexity} $c(M,X)$ as the number of vertices of any minimal
skeleton of $(M,X)$.

\paragraph{Examples with complexity zero.}

\begin{itemize}

\item It is well-known~\cite{Mat} that the only closed prime manifolds
having complexity zero are $S^3, S^2\times S^1, \matRP^3$, and $L_{3,1}$.

\item The trivial element $B_0=(T\times [0,1],\{\theta \times \{0\},
\theta \times \{1\}\})$ has complexity zero, since it has the simple
skeleton $\theta\times [0,1]\subset T\times [0,1]$, which has no vertices.

\item Let $H$ be the solid torus, let $D$ be a meridinal disc properly
embedded in $H$ and let $\theta \subset\partial H$ be a triod containing
$\partial D$, as in Fig.~\ref{smallbr}-left. Then $D \cup \theta$ is a
skeleton of $B_1=(H,\{\theta\})$, which has therefore complexity zero.

\begin{figure} \begin{center}
\includegraphics{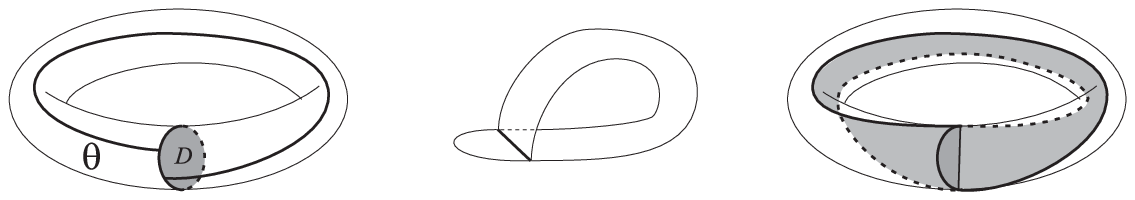}
\nota{The minimal skeleta of $B_1$ and $B_2$.} \label{smallbr}
\end{center} \end{figure}

\item Let $H$ be the solid torus again, and let $P$ be the M\"obius strip
with one tongue shown in Fig.~\ref{smallbr}-centre, embedded in $H$ as in
Fig.~\ref{smallbr}-right. Since $P$ has no vertices and it is a skeleton
for $B_2=(H,\{P\cap\partial H\})$, then $c(B_2)=0$.

\end{itemize}

\subsection{Properties of complexity}\label{ct-properties} Of course we
have $c(M,\emptyset)=c(M)$, namely property~(\ref{extension:property}) of
our list. We prove in this subsection the other properties of $c$. This
will require, together with some \emph{ad hoc} methods, the extension to
our context of some techniques used in~\cite{Mat}. In the course of our
arguments we will give several definitions used elsewhere in the paper,
and we will prove other facts stated above.

\paragraph{Finiteness.} The proof of property~(\ref{finite:property}) of
complexity requires a careful discussion of the topological properties of
minimal spines.

\begin{figure}
\begin{center}
\includegraphics{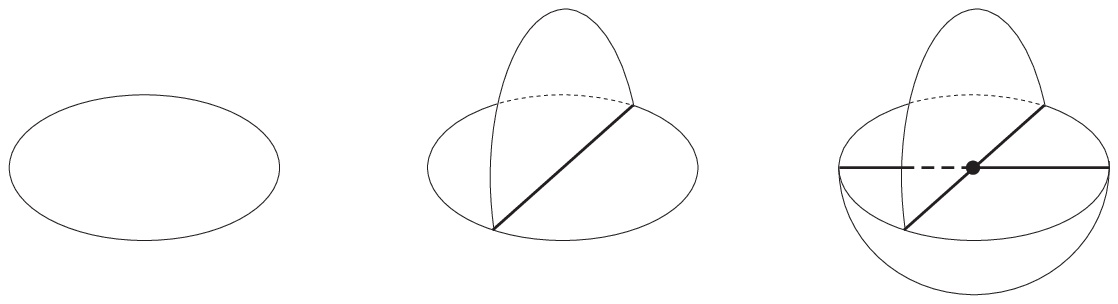}
\nota{Typical neighbourhoods of points in a quasi-standard
polyhedron.} \label{neighb}
\end{center}
\end{figure}

A simple polyhedron $Q$ is called \emph{quasi-standard} if the link
of every point is either a circle, or a circle with a diameter, or a
circle with three radii (neighbourhoods of points of the three types
are shown in Fig.~\ref{neighb}).
A simple polyhedron $Q$ is called \emph{quasi-standard with boundary} if
in addition to these three types of points we have points having as a link
either a closed segment or the union of $3$ closed segments with one
common endpoint. Assuming $Q$ to be quasi-standard with boundary, we
denote by $V(Q)$ the set of points (called \emph{vertices} above) whose
link is a circle with three radii, and by $S(Q)$ the union of $V(Q)$ with
the set of points whose link is a circle with a diameter. We also denote
by $\partial Q$ the points of the two new types declared legal when
passing from `quasi-standard' to `quasi-standard with boundary.' Moreover,
we call \emph{$1$-components} of $Q$ the connected components of $S(Q)
\setminus V(Q)$ and \emph{$2$-components} of $Q$ the connected components
of $Q \setminus (S(Q)\cup\partial Q)$. If the $2$-components of $Q$ are open discs (and
hence are called just \emph{faces}), and the $1$-components are open
segments (and hence called just {\em edges}), then we call $Q$ a
\emph{standard polyhedron with boundary}. For short we will often just
call $Q$ a \emph{standard} polyhedron, and possibly specify that $\partial
Q$ should or not be empty.

We state now several easy facts concerning nuclear skeleta, and prove
a crucial result concerning minimal skeleta.

\begin{rem} \label{nuclear:rem} \emph{Let $(M,X)$ be a manifold with
triods and let $P$ be a nuclear skeleton of $(M,X)$. Then, up to
rearranging the components $T_1,\ldots,T_n$ of $\partial M$, we have that
$P=Q\cup s_1\cup\ldots\cup s_m\cup K$, where:
\begin{enumerate}
\item $Q$ is a quasi-standard polyhedron with boundary $\partial Q \subset X$;
\item For $i=1,\ldots,m$ we have that $s_i\subset\theta_i$ is a segment
and $Q\cup s_i$ appears near $T_i$ precisely as the minimal skeleton of
$B_1$ appears near $\partial B_1$ (see Fig.~\ref{smallbr}-left); for $i>m$
we have $\partial Q\supset\theta_i$;
\item $K$ is a graph with $K\cap (Q\cup s_1\cup\ldots\cup s_m)$ finite and
$K\cap V(Q\cup\partial M)$ empty.
\end{enumerate}} \end{rem}

\begin{rem}\label{nuclear:rem:bis} \emph{Every $(M,X)\in\calX$ has a
minimal skeleton $P'=Q\cup s_1\cup\ldots\cup s_m\cup K'$ as above,
where in addition $K'\cap\partial M=\emptyset$. This is because,
without changing $\#V(P)$, we can take the ends of $K$ lying on
$\partial M$ and make them slide over $Q\cup s_1\cup\ldots\cup s_m$
until they reach $\interior{M}$. Note that the regular neighbourhood
of $\theta_i \in X$ in $P'$ is now either a product $\theta_i\times
[0,1]$ or the union of an annulus and a segment, as for
$B_1$.}\end{rem}

\begin{rem}\label{standard:is:proper} \emph{If $P$ is a nuclear and
standard skeleton of $(M,X)$ then it is properly embedded, namely
$\partial P =\partial M \cap P = X$, and $P\cup\partial M$ is standard
without boundary. Moreover $P\cup\partial M$
is a spine of a manifold bounded by one sphere and some tori, so
$\chi(P\cup\partial M)=1$. Knowing that $S(P\cup\partial M)$ is 4-valent
and denoting by $F(P)$ the set of faces of $P$, we also see that
$\#F(P)-\#V(P)=\#X+1$.}\end{rem}

\begin{teo} \label{standard:teo} Let $(M,X) \in \calX$ be prime and
let $P$ be a minimal skeleton of $(M,X)$. Then:

\begin{enumerate}
\item If $c(M,X)>0$ then $P$ is standard;
\item If $c(M,X)=0$ and $X\neq\emptyset$ then $(M,X)\in\{B_0, B_1, B_2\}$, and
$P$ is the skeleton described in Subsection~\ref{ct-def} (which is
standard for $B_0$ and $B_2$ only);
\item If $c(M,X)=0$ and $X=\emptyset$ then $(M,X)\in\{S^3,S^2\times
S^1,L_{3,1},\matRP^3\}$ and $P$ is not standard.
\end{enumerate}
\end{teo}

\begin{proof} Our argument closely follows~\cite{Mat}.
We can first rule out the case $(M,X)=(S^2\times S^1,\emptyset)$,
because for it we only need to show that $P$ is not standard. But a
standard polyhedron without boundary must have vertices, while
$c(S^2\times S^1,\emptyset)=0$. So we proceed assuming that $M$ is
irreducible.

We will now prove that if $P$ is not standard then
$(M,X)\in\{B_1,S^3,L_{3,1},\matRP^3\}$, and that $P$ is as in
Subsection~\ref{ct-def} when $(M,X)=B_1$. To conclude we will later
show that if $P$ is standard and $c(M,X)=0$ then
$(M,X)\in\{B_0,B_2\}$ and $P$ is as prescribed.

Suppose then $P$ is not standard. First, if $P$ is a point then
$(M,X)=(S^3,\emptyset)$. Suppose now $P$ has a 1-dimensional part.
So, let
 $e\subset P$ be a segment disjoint from the 2-dimensional part of
$P$. If $e\subset\partial M$, looking at the ball
$M\setminus(P\cup\partial M)$, we deduce that there is a properly
embedded disc in $M$ intersecting $P$ in a point of $e$. By
irreducibility $M$ is then a solid torus, so $(M,X)=B_1$ and $P$ is
as in Fig.~\ref{smallbr}-left. If $e\subset\interior M$, looking at
the ball $M\setminus(P\cup\partial M)$ again, we see that there is a
sphere $S\subset M$ intersecting $P$ in one point of $e$. By
irreducibility $S$ bounds a ball $B$, and $P\cap B$ is easily seen to
be a spine of $B$. Nuclearity now implies that $P\cap B$ contains
vertices, so $P\setminus B$ is a skeleton of $(M,X)$ with fewer
vertices than $P$. A contradiction.

We have shown so far that $P$ is quasi-standard unless $(M,X)$ is
$S^3$ or $B_1$. Since $P$ is not standard, either a 2-component $f$ is not
a disc, or a 1-component is a circle. In the first case, either $f=S^2$, or
$f=\matRP^2$, or $f$ contains a simple closed curve $\gamma$ which is
non-trivial and orientation-preserving in $f$. In the first two cases
we have respectively $P=S^2$, which is impossible, and $P=\matRP^2$,
so $(M,X)=\matRP^3$. The third case is impossible: looking once more
at the ball $M\setminus(P\cup\partial M)$, we deduce that there is a
sphere $S\subset M$ intersecting $P$ in $\gamma$, and again
$S=\partial B$. As above, $P\cap B$ is a spine of $B$. By minimality
$P\cap B$ cannot contain vertices. It follows that $P\cap B$ is a
disc, which contradicts the choice of $\gamma$.
Finally, if a 1-component of $P$ is a circle but all 2-components are
discs, then $P$ must be the ``triple hat,'' a skeleton of $L_{3,1}$.

We are left to analyze the case where $P$
is standard and $c(M,X)=0$, so $X\ne\emptyset$.
Now, if $\theta\in X$ and $p$ is a vertex
of $\theta$, then the three faces of $P$ incident to $p$ are the same
as those incident to the other vertex of $\theta$. Moreover, since
$V(P)=\emptyset$, again the same faces are incident to the endpoint
of the edge of $P$ which starts at $p$. It easily follows that
$F(P)\le 3$, but $F(P)=1+\#X$ by Remark~\ref{standard:is:proper}, so
$\#X$ is either $1$ or $2$. It is now a routine matter to check that
$(M,X)$ is respectively $B_2$ or $B_0$, with $P$ as
prescribed.\end{proof}

The next two results show respectively
property~(\ref{finite:property}) of complexity and
Lemma~\ref{zero:compl:lem}.

\begin{cor} For any $n\ge 0$, only finitely many pairs in $\calX^{\rm pr}$
have complexity $n$. \end{cor}

\begin{cor}\label{zero:bricks} $\calB_0^0=\emptyset$ and
$\calB_0=\calB_0^1=\{B_0,B_1,B_2\}$. \end{cor}

\begin{proof} There are no closed bricks of complexity
zero, since $(S^2\times S^1, \emptyset)$, $(S^3,\emptyset)$,
$(\matRP^3,\emptyset)$, and $(L_{3,1},\emptyset)$ can be obtained
assembling respectively two copies of $B_1$, two copies of $B_1$, one copy of
$B_1$ and one of $B_2$, and two copies of $B_2$. Moreover $B_0$, $B_1$, and $B_2$
are not non-trivial
assemblings of each other, and the conclusion follows.\end{proof}

\paragraph{Subadditivity under (self-)assembling.} Let $(M,X)$ and
$(M',X')$ be two given pairs, and let $(N,Y)$ be obtained by assembling
them. Let $P$ and $P'$ be minimal skeleta respectively of $(M,X)$ and
$(M',X')$. The assembling is defined by an identification $\psi:T_i \to
T_{i'}'$ with $\psi(\theta_i) = \theta_{i'}'$. Using
Remark~\ref{nuclear:rem} we see that $P\cup_{\psi} P'$ is simple, so it is
a skeleton of $(N,Y)$, and that no new vertices appear. It follows that
$c(N,Y)\le c(M,X)+c(M',X')$.

Let $(M,X)$ be a pair and let $(N,Y)$ be obtained from $(M,X)$ via a
self-assembling, determined by a map $\psi:T_i \to T_{i'}$ such that
$\psi(\theta_i)$ intersects transversely $\theta_{i'}$ in two points. If
$P$ is a minimal skeleton of $(M,X)$ as in Remark~\ref{nuclear:rem:bis},
then $P\cup T_i \subset N$ is a skeleton for $(N,X)$. Moreover $P\cup T_i$
has at most 6 vertices more than $P$ (2 from the vertices of $\theta_i$, 2
from those of $\theta_{i'}$, and 2 from $\psi(\theta_i)\cap\theta_{i'}$).
It follows that $c(N,Y)\le c(M,X)+6$.

\paragraph{Normal surfaces.} Let $(M,X)$ be a manifold with triods and let
$P$ be a nuclear skeleton of $(M,X)$. The simple polyhedron
$P\cup\partial M$ is now a spine of $M$ with a ball $B\subset M$
removed. Choose a triangulation of $P\cup\partial M$, and let $\xi_P$
be the handle decomposition of $M\setminus B$ obtained thickening the
triangulation of $P\cup\partial M$, as in~\cite{Mat}. In this
paragraph we will study closed normal surfaces in $\xi_P$. A
connected normal surface $S$ is \emph{parallel to the boundary} when
it is obtained by taking one boundary component and pushing it a bit
inside $\xi_P$. In our case, we have one such surface for each $T_i$,
and one for $\partial B$.

Two preliminary results are needed to prove our main statement on
normal surfaces. The first one refers to another situation, very
often considered below, where a normal surface naturally arises.

\begin{prop} \label{loops-prop} Let $(M,X)$ be a manifold with triods and
let $Q\subset M$ be a quasi-standard polyhedron with
$Q\cap\partial M=\partial Q\subset
 X$.
Assume $M\setminus Q$ has two
components $N'$ and $N''$. Then the
 faces of $Q$ that separate $N'$ from
$N''$ form a closed orientable
 surface $\Sigma(Q) \subset Q \subset M$ which
cuts $M$ into two
 components. \end{prop}

\begin{proof} Let $e$ be an edge of $Q$, and let $\{f_1, f_2, f_3\}$ be
the triple of (possibly not distinct) faces of $Q$ incident to $e$. The
number of $f_i$'s that separate $N'$ from $N''$ is even; it follows that
$\Sigma(Q)$ is a surface away from $V(Q)\cup\partial Q$. Let $T_i$ be a
boundary component of $M$, containing the triod $\theta_i \in X$. Since
$T_i\setminus \theta_i$ is a disc, which is adjacent either to $N'$ or to
$N''$ (say $N'$), then each $2$-component of $Q$ incident to $\theta_i$
(there could be $0$, $1$ or 3 of them, with multiplicity) has $N'$ on both
sides. So $\Sigma(Q)$ is not adjacent to $\partial Q$. Finally, since
$\Sigma(Q)$ intersects the link of each vertex either nowhere or in a
loop, then $\Sigma(Q)$ is a closed surface.

The surface $\Sigma(Q)$ cuts $M$ in two components (and is thus
orientable, since $M$ is) because $N'$ and $N''$ lie on opposite sides of
$\Sigma(Q)$. \end{proof}

\begin{lemma} \label{standard:has:vertices} Let $P$ be a standard and
nuclear skeleton of a pair $(M,X)$. If $\#V(P)>0$ then every face of
$P$ is incident to at least one vertex. \end{lemma}

\begin{proof} Assume a face $f$ of $P$ contains no vertices, and let $f$
be incident to the triods $\theta_{i_1},\ldots,\theta_{i_k}$. Then
$\partial f\cup\theta_{i_1}\cup\ldots\cup\theta_{i_k}$ is a connected
component of $S(P\cup\partial M)$, but $P\cup\partial M$ is standard
without boundary by Remark~\ref{standard:is:proper}, so
$S(P\cup\partial M)=\partial
f\cup\theta_{i_1}\cup\ldots\cup\theta_{i_k}$, whence
$S(P)\subset\partial f$ and $V(P)=\emptyset$. A contradiction.
\end{proof}

We go back now to the situation where $P$ is a nuclear skeleton of
$(M,X)$.

\begin{lemma} \label{normal:lemma} Let $F$ be a closed normal surface in
$\xi_P$. Assume that no component of $F$ is boundary-parallel. Then
there exists a simple polyhedron $P_F$ embedded in $M$, with
$\#V(P_F)\le \#V(P)$, such that $P_F\cap \partial M = X$ and
$M\setminus (P_F\cup\partial M)$ is an open regular neighbourhood of
$F$. Moreover, if $P$ is standard and $\#V(P)>0$ then $\#V(P_F) <
\#V(P)$. \end{lemma}

\begin{proof} Being normal, $F$ is determined by an integer attached to
each 2-component of $P\cup\partial M$. Now we cut $P\cup\partial M$
open along $F$ as explained in~\cite{Mat}: if a 2-component bears an
integer $n$ we replace the component by $n+1$ parallel ones. We get a
polyhedron $P'\subset M$ which contains $\partial M$, such that
$M\setminus P'$ is the disjoint union of an open ball $B$ and an open
regular neighbourhood $N$ of $F$ in $M$. By removing from each torus
$T_i\subset \partial M$ the open disc $T_i\setminus\theta_i$ we get a
polyhedron $P''$ intersecting $\partial M$ in $X$. Now we puncture a
$2$-component which separates $B$ from $N$ and claim that the
polyhedron $P_F$ is as desired. Only the inequalities between $V(P)$
and $V(P_F)$ are non-obvious.

By construction we have $\#V(P_F\cup\partial M) \le \#V(P\cup\partial M)$.
Consider now a vertex $v$ of $P\cup\partial M$ contained in
$T_i\subset\partial M$. Of the six germs of 2-component of $P\cup\partial
M$ at $v$, three are actually the same $T_i\setminus\theta_i$, so their
coefficient in $F$ is the same, say $\alpha$. Call $\beta$, $\gamma$, and
$\delta$ the coefficients of the other three germs of 2-component at $v$.
As we cut $P\cup\partial M$ along $F$ we see that $v$ disappears if and
only if (up to permutation) $\beta=\gamma>\delta$. If $v$ does not
disappear then $\beta=\gamma=\delta$ is even. Then we set
$k=\alpha-\beta/2$ and note that $v$ remains on $\partial M$ if and only
if $k=0$. Now let $v'$ be the other vertex of $P\cup\partial M$ on $T_i$.
Since the coefficients $(\alpha,\alpha,\alpha,\beta,\gamma,\delta)$ are
the same at $v'$, we deduce that either $v$ and $v'$ both disappear, or
they both stay on $\partial M$, or they both move to $\interior{M}$. In
the last case, however, one sees that $F$ has $k$ components parallel to
$T_i$, which is absurd. So both $v$ and $v'$ disappear in $P''$ (either
already in $P'$ or when we remove $T_i\setminus\theta_i$). This shows that
$\#V(P'')\le\#V(P)$, so $\#V(P_F)\le\#V(P)$.

Suppose now $P$ is standard. Then $P''$ is the union of a quasi-standard
polyhedron $P'''$ and some arcs in $X$.
The 2-components of $P''$ which separate $B$ from $N$ are the same as those
of $P'''$, so they give a closed surface $\Sigma\subset P''$ by
Proposition~\ref{loops-prop}.
Since no component of $F$ is parallel to $\partial B$
or to one of the $T_i$'s, the 2-component $f$ of $P''$ punctured to
get $P_F$ cannot be a closed surface. Now if $\partial f$ contains
vertices of $P''$, we see that  $\#V(P_F)<\#V(P'')\le\#V(P)$, whence
the conclusion. Suppose on the
contrary that $\partial f$ contains a circle $\gamma\subset S(P'')$ with
$\gamma\cap V(P'')=\emptyset$. Note that the process of cutting $P\cup\partial
M$ along $F$ allows to define a local injection $\psi:P'\to P\cup\partial
M$, and that $P''\subset P'$. Now, if $\psi(\gamma)$ contains some vertex
of $P$ then this vertex has disappeared in the passage from $P$ to $P''$,
whence the conclusion. If $\psi(\gamma)\cap V(P)=\emptyset$ then we
consider the 2-component $g$ of $P''\setminus\Sigma$ incident to
$\gamma$ and note that $\psi(g)$ must be a face of $P$ without vertices,
which is absurd by Lemma~\ref{standard:has:vertices}.\end{proof}

\begin{teo} \label{standard2:teo} If $(M,X)\in\calX$  has a
standard minimal skeleton then it is prime. \end{teo}

\begin{proof} For $c(M,X)=0$ it was shown during the proof of
Theorem~\ref{standard:teo} that $(M,X)$ is $B_0$ or $B_2$, so we suppose
$c(M,X)>0$. By contradiction, assume $M$ is not prime and let $P$ be a
standard minimal skeleton of $(M,X)$. Then $\xi_P$ contains an essential
normal sphere $S$. Such a sphere cannot be parallel to the boundary in
$\xi_P$. Applying Lemma~\ref{normal:lemma} we get $P_S\subset M$ with
$\#V(P_S)<\#V(P)$, $P_S\cap\partial M=X$, and $M\setminus (P_S\cup\partial
M) \cong S\times (0,1)$. Since $(S\setminus\{{\rm point}\})\times(0,1)$ is
an open 3-ball, adding to $P_S$ a generic segment isotopic to $\{{\rm
point}\}\times(0,1)$ we get a skeleton for $(M,X)$ with as many vertices
as $P_S$. This contradicts minimality of $P$. \end{proof}

\paragraph{Additivity under connected sum.} Again, we follow \cite{Mat}
quite closely. Let $(M,X)$ and $(M',X')$ be manifolds with triods, and set
$(N,Y)=(M,X)\#(M',X')$. Let $P$ and $P'$ be skeleta of $(M,X)$ and
$(M'X')$, respectively. If we take points $p\in P$ and $p'\in P'$ which
are not vertices and we join them with a segment, we get a skeleton of
$(N,Y)$. This implies that $c(N,Y)\le c(M,X)+c(M',X')$.

Let us prove the opposite inequality. Let $P$ be a minimal skeleton of
$(N,Y)$. Since $(N,Y)$ is not prime, there is a separating normal sphere
$S$ in $\xi_P$ (maybe not the one which cuts $N$ into $M$ and $M'$, as
customary in normal surface theory). Let $(N_1,Y_1)$ and $(N_2,Y_2)$ be
obtained by cutting $(N,Y)$ along $S$ and gluing in balls. The polyhedron
$P_S$ given by Lemma~\ref{normal:lemma} is now the disjoint union of two
polyhedra $P_1$ and $P_2$ such that $P_i$ is a skeleton of $(N_i, Y_i)$.
Moreover $\#V(P_S)=\#V(P_1)+\#V(P_2)\le \#V(P)$. Therefore
$c(N_1,Y_1)+c(N_2,Y_2)\le c(N,Y)$, whence $c(N_1,Y_1)+c(N_2,Y_2)=c(N,Y)$.
We can now go on finding essential spheres, and additivity eventually
follows from uniqueness of the decomposition into primes.

\paragraph{Sharp (self-)assemblings.} We are now in a position to prove
the second half of properties~(\ref{assemble:subad:property})
and~(\ref{self-assemble:subad:property}) of complexity. The case of
self-assembling is actually easier, so we start from it. Let a sharp
$(N,Y)=\odot(M,X)$ be performed along $\psi:T_i\to T_{i'}$.
Let $P$ be a minimal skeleton of $(M,X)$ as in
Remark~\ref{nuclear:rem:bis}. Then $P\cup T_i$ is a minimal skeleton
of $(N,Y)$, and it is easy to see that $P$ is standard if and only if
$P\cup T_i$ is. Moreover, by Theorem~\ref{standard:teo} and
Theorem~\ref{standard2:teo}, $P$ is standard if and only if $(M,X)$
is prime (because $\# X\ge 2$) and $P\cup T_i$ is standard if and
only if $(N,Y)$ is prime (because $c(N,Y)>0$).
This shows the desired conclusion that $(M,X)$ is prime if
and only if $(N,Y)$ is.

To deal with assembling, we need two preliminary results. The first
one, together with Theorem~\ref{standard:teo}, implies
Proposition~\ref{not:sharp}.

\begin{lemma}\label{B1:not:used} Let $(M,X)\in\calX$ be prime and assume
$c(M,X)>0$. Then no assembling $(M,X)\oplus B_1$ is sharp.\end{lemma}

\begin{proof} Let $P$ be a
minimal skeleton for $(M,X)$, which is standard by
Theorem~\ref{standard:teo}, and let $P'$ be the minimal skeleton of $B_1$.
Then $P\cup_{\psi}P'$ is a skeleton for $(M,X)\oplus B_1$ with minimal number of
vertices, but $P\cup_{\psi}P'$ is not nuclear: there is a face $f$ of $P$,
glued to the free segment of $P'$, which is incident to some vertex of $P$ by
Lemma~\ref{standard:has:vertices}. By collapsing $f$ we would get a
skeleton with fewer vertices, which is absurd. \end{proof}

\begin{lemma} \label{faces:distinct} Let $P$ be a minimal skeleton of
$(M,X)\in\calX^{\rm pr}$ with $c(M,X)>0$. Then, for each
$\theta_i\in X$, the three faces of $P$ incident to $\theta_i$ are
distinct from each other. \end{lemma}

\begin{proof} By Theorem~\ref{standard:teo}, $P$ is standard. Suppose a
face $f$ is incident more than once to some $\theta_i$. Let $\alpha$
be an arc in $f$ having endpoints $p_0$ and $p_1$ in two distinct
edges of $\theta_i$, and let $\beta$ be an essential closed curve in
$T_i\subset\partial M$ with $\beta\cap\theta_i = \{p_0,p_1\}$. Now
$\beta$ is cut by $\{p_0,p_1\}$ into components $\beta_0$ and
$\beta_1$. Since $M\setminus(P\cup\partial M)$ is a ball, we can glue
to both curves $\alpha \cup \beta_i$ a disc, and the two discs
together form a disc $D\subset M$ with $\partial D = \beta$. Since
$\beta$ is essential, $M$ is a solid torus and $(M,X)=B_1$.
\end{proof}

Now let $(N,Y)=(M,X)\oplus(M',X')$ be a sharp assembling along some
map $\psi:T_i\to T'_{i'}$. Recall that we want to show that $(N,Y)$
is prime if and only if both $(M,X)$ and $(M',X')$ are. Assume first
that $c(N,Y)=0$. If $(M,X)$ and $(M',X')$ are prime, by
Theorem~\ref{standard:teo} $(N,Y)$ is a lens space, so it is prime.
If $(N,Y)$ is prime, we consider the prime factorization of $(M,X)$
and $(M',X')$, and note that $\psi$ assembles one factor $W$ of
$(M,X)$ to one factor $W'$ of $(M',X')$. If $W\oplus W'\ne
(S^3,\emptyset)$, then, since $(N,Y)$ is prime, $(M,X)=W$ and
$(M',X')=W'$, and we are done. Otherwise, up to permutation,
$(M,X)=Z\# W$ and $(M',X')=W'$. By additivity of $c$ under $\#$ and
Theorem~\ref{standard:teo}, $W$ and $W'$ are solid tori, and the
assembling is trivial.

Now let $c(N,Y)$ be positive. Up to permutation, $c(M,X)>0$.
Let $P$ and $P'$ be minimal skeleta of $(M,X)$ and $(M',X')$
respectively, so $P\cup_{\psi} P'$ is a minimal skeleton of $(N,Y)$.
If $(N,Y)$ is prime, $P\cup_{\psi} P'$ is standard by
Theorem~\ref{standard:teo}, so
$P$ and $P'$ are, and Theorem~\ref{standard2:teo}
implies the conclusion. Conversely, let $(M,X)$ and $(M',X')$ be prime.
If $(M',X')=B_1$, we get a contradiction to Lemma~\ref{B1:not:used}.
Otherwise Theorem~\ref{standard:teo} implies that $P$ and $P'$ are standard.
Now, it is not a priori obvious that $P\cup_{\psi} P'$ is standard, because
some annular component could appear, but
Lemma~\ref{faces:distinct} applied to $P$ shows that they actually do not,
and our argument is complete.

\begin{rem}\label{non-trivial:rough}
\emph{Given $(H,\{\theta\})\in\calX$ with $H$ the solid torus, it is
easy to see that there are infinitely many $(H,\{\theta'\})$'s such
that $(H,\{\theta\})\oplus(H,\{\theta'\})=(S^3,\emptyset)$, so
$((M,X)\#(H,\{\theta\}))\oplus(H,\{\theta'\})=(M,X)$ for any $(M,X)$.
However, the only assemblings of this sort on which complexity is
additive are those where $c(H,\{\theta\})=c(H,\{\theta'\})=0$. This
can only happen if $\{(H,\{\theta\}),(H,\{\theta'\})\}$ is
$\{B_1,B_1\}$ or $\{B_1,B_2\}$, so these are the only cases which our
definition of \emph{trivial} rules out from the notion of
\emph{sharp} assembling.}\end{rem}

\section{The algorithm to find bricks} \label{algorithm} We will explain
in this section how we have been able to determine $\calB_{\le 9}$.

\subsection{Properties of minimal skeleta of bricks}
We will introduce in this subsection
two more bricks $B_3$ and $B_4$, besides the $B_0$,
$B_1$ and $B_2$ already defined above. Then we will state some results
giving strong restrictions on the shape of minimal skeleta of bricks
different from $B_0,\ldots,B_4$. Later we will describe the
operations which we actually have carried out by computer to
determine $\calB_{\le 9}$.

\paragraph{Minimal skeleta for $B_3$ and $B_4$.} We define $B_3$ and $B_4$
as the elements of $\calX$ based on $D_1 \times S^1$ and $D_2 \times
S^1$ respectively, where $D_i$ is the disc with $i$ holes, and the
boundary triods are as decribed in Table~\ref{brick:table}
(Subsection~\ref{results}). A skeleton for $B_3$ is given by the
union of an annulus $D_1\times\{{\rm point}\}$ and a ribbon, glued as
in Fig.~\ref{smallbr2}-left. Similarly, a skeleton for $B_4$ is given
by the union of $D_2\times\{{\rm point}\}$ and a polyhedron as in
Fig.~\ref{smallbr2}-right, glued as shown. This implies that
$c(B_3)\le 1$ and $c(B_4)\le 3$. Since $B_3$ is prime and it is not
$B_0$, $B_1$, or $B_2$, we have $c(B_3)=1$ by
Theorem~\ref{standard:teo}-(2). Using Theorem~\ref{standard:teo}-(1)
and checking by hand all standard $P$'s with $\#V(P)=1$ and $\partial
P\ne\emptyset$, we see that $B_3$ is a brick and actually
$\calB_1^1=\{B_3\}$. For $B_4$ we need:

\begin{figure}
\begin{center}
\includegraphics{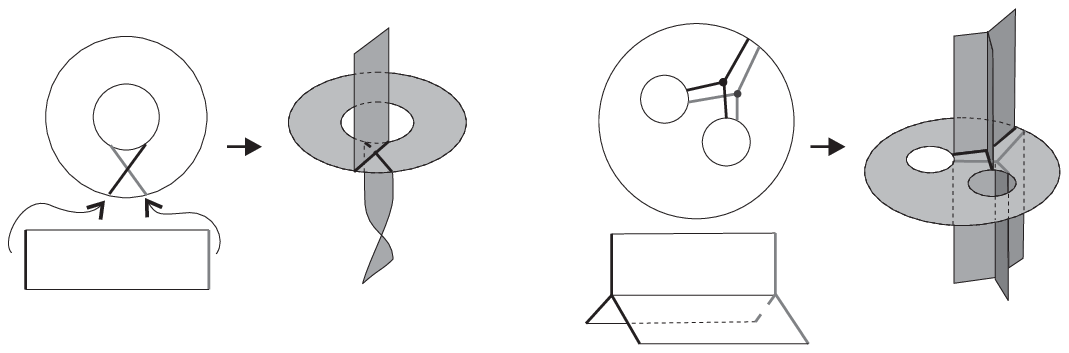}
\nota{Minimal skeleta for $B_3$ and $B_4$. In the 3-dimensional (gray)
pictures the segment and the $Y$ on the top are
identified with the corresponding segment and $Y$ on the bottom.} \label{smallbr2}
\end{center} \end{figure}

\begin{lemma} Let $(M,X)$ be prime and different from $B_0,\ldots,B_3$. Then
$c(M,X)\ge \#X$. \end{lemma}

\begin{proof}
Of course we can assume $X\ne\emptyset$.
Since $\calB_{\le1}^1=\{B_0,\ldots,B_3\}$ and the inequality is easy for any
non-trivial assembling of $B_0,\ldots,B_3$, we also assume $c(M,X)\ge2$.

Suppose now that a face $f$ is incident to two distinct triods
$\theta_i, \theta_{i'} \in X$. Then there is an arc $\lambda\subset
f$, properly embedded in $f$, with endpoints $p\in\theta_i,
p'\in\theta_{i'}$, and two essential loops $\gamma\subset T_i$ and
$\gamma'\subset T_{i'}$ such that $\gamma\cap\theta_{i}=\{p\}$,
$\gamma'\cap\theta_{i'}=\{p'\}$. Since
 $M\setminus(P\cup\partial M)$ is a
ball, there is an annulus
 $A$ properly embedded in $M$,
with $\partial A=\gamma\cup\gamma'$ and $A\cap P=\lambda$. If some
face $g\neq f$ is incident to the same $\theta_i$ and to some other
$\theta_{i''}$, we can construct an annulus $B$ in the same way.
Moreover $\partial B=\delta\cup\delta''$ with
$\#(\gamma\pitchfork\delta)=1$. Irreducibility allows to assume that
$A\cap B$ is just one segment, hence $\theta_{i'}=\theta_{i''}$, and
then to show that $M=T\times[0,1]$. So $\#X=2$, but we are assuming
$c(M,X)\ge2$, and the conclusion holds in this case.

By Lemma~\ref{faces:distinct}, $P$ has distinct faces
$f_i^{(1)},f_i^{(2)},f_i^{(3)}$ incident to $\theta_i$. By what
already shown we can assume up to permutation that $f_i^{(2)}$ and
$f_i^{(3)}$ are not incident to any other triod in $X$. So $P$
contains at least $2\#X+1$ distinct faces. By
Remark~\ref{standard:is:proper} we have $\#F(P)=\#V(P)+\#X+1$,
therefore $\#V(P)\ge \#X$.
\end{proof}

\begin{prop}
$c(B_4)=3$ and $B_4$ is a brick.
\end{prop}

\begin{proof}
First, $B_4$ is prime and $\#X(B_4)=3$, so $c(B_4)=3$ by the previous
lemma. If $B_4$ were not a brick then it would split as $B_i\oplus
B^{(1)}\oplus\ldots\oplus B^{(k)}$ with $i\in\{1,2,3\}$. In all cases
we must have $\#X(B^{(j)})>c(B^{(j)})$ for some $j$, which
contradicts the previous lemma.\end{proof}

\paragraph{Super-standard skeleta.}

A standard polyhedron $P$ (with boundary) is called \emph{super-standard} if
every face of $P$ is incident to $\partial P$ along one segment at most.
For such a $P$, it is easy to prove that $S(P)$ must be
connected if $P$ is. The minimal skeleta of $B_0,\ldots,B_4$ we have
already described are not super-standard. The following theorem will be
proved in Section~\ref{proofs:section}.

\begin{teo} \label{minimal:are:standard} Let $(M,X)$ be a brick different
from $B_0,\ldots,B_4$. Then every minimal skeleton of $(M,X)$ is
super-standard. \end{teo}

\begin{cor} Let $(M,X)$ be a brick different from $B_0,\ldots,B_4$. Then
$c(M,X)\ge 2\#X-1$. \end{cor}

\begin{proof} Let $P$ be a minimal skeleton of $(M,X)$.
By Remark~\ref{standard:is:proper} we have $\#F(P)-\#V(P)=\#X+1$. Now
$P$ is super-standard, so each edge in $X$ determines a different
face of $P$. Then $\#F(P) \ge 3\#X$, and the conclusion follows.
\end{proof}

\paragraph{Enumeration of bricks.} Let $(M,X)$ be a brick different from
$B_0,\ldots,B_4$, and let $P$ be one of its minimal skeleta. We will call
\emph{filling} of $P$ any of the (finitely many) polyhedra obtained by
glueing to $P$ one copy of the M\"obius strip with one tongue along each
of the boundary triods of $P$ (so $\#X$ strips in all are glued to $P$).
Since the M\"obius strip with one tongue is a skeleton of the pair $B_2$
based on the solid torus,
we see that a filling of $P$ is automatically a skeleton of a
(possibly non-sharp)
assembling of $(M,X)$ with $\#X$ copies of $B_2$, hence of a closed
manifold $(N,\emptyset)\in\calX$ obtained by Dehn-filling all boundary
components of $M$. Note that the glueing function $\psi$ used to define
the filling of one component $T_i$ of $\partial M$ must map the triod
$\theta_i\subset T_i$ to the triod of $B_2$, so indeed there are finitely
many possibilities. Since $P$ is super-standard by
Theorem~\ref{minimal:are:standard}, it is easy to see that the fillings of
$P$ are standard.

We will call \emph{loop} in $P$ a subpolyhedron $\gamma \subset P$
homeomorphic to $S^1$ which intersects $S(P)$ transversely (in
particular $\gamma \cap V(P) = \emptyset)$. We define the
\emph{length} $l(\gamma)$ of $\gamma$ as the number of its
intersections with the edges of $P$. We denote by $\calR(\gamma)$ a
regular neighbourhood of $\gamma$ in $P$. Note that the core of the
M\"obius strip has length 1 in the M\"obius strip with one tongue.
The following result will be shown in Section~\ref{proofs:section}.

\begin{teo} \label{disjoint:length:1:loops} Let $(M,X)$ be a brick
different from $B_0,\ldots,B_4$. Let $P$ be a minimal skeleton of $(M,X)$
and let $Q$ be any filling of $P$. Let $\calL(Q)$ be any set of
representatives of the ambient isotopy classes of length-1 loops in $Q$.
Then:
\begin{enumerate}
\item The elements of $\calL(Q)$ are pairwise disjoint, and $\calR(\gamma)$
is a M\"obius strip with one tongue for all $\gamma\in\calL(Q)$;
\item $\calL(Q)$
consists of $\#X$ loops and $P = Q\setminus \calR(\calL(Q))$.
\end{enumerate}
\end{teo}

\begin{rem}\label{real:meaning:loops}
\emph{Condition (1) in Theorem~\ref{disjoint:length:1:loops} means that:
\begin{itemize}
\item for every edge $e$ of $Q$, there is no face $f$ of $Q$ triply incident
to $e$;
\item if $f$ is doubly incident to $e$ and $\partial f$ is oriented, then
$e$ is induced the same orientation twice;
\item for $i=1,2$ let $f_i$ be doubly incident and $g_i$ be incident to
the same $e_i$, with $e_1\ne e_2$; then $f_1\ne f_2$.
\end{itemize}
In addition, taking point (2) of Theorem~\ref{disjoint:length:1:loops}
for granted, super-standardness of $P$ (stated by
Theorem~\ref{minimal:are:standard}) means the following:
\begin{itemize}
\item with the above notation, $f_1,f_2,g_1,g_2$ are pairwise distinct.
\end{itemize}}
\end{rem}

We state now a result on the singular set $S(Q)$ of a filling $Q$ of a
minimal skeleton $P$, noting first that $S(Q)$ depends on $P$ only and it
is a 4-valent graph (because $\partial Q=\emptyset$). We refer again to
Section~\ref{proofs:section} for the proof.

\begin{figure} \begin{center}
\includegraphics{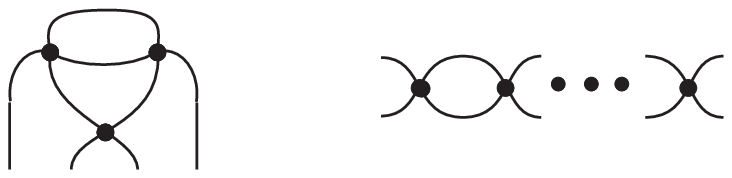}
\nota{If 4 edges disconnect $S(Q)$, then one of the two pieces
is of one of these types.} \label{2bridge}
\end{center} \end{figure}

\begin{teo} \label{restrictions:on:graph} Let $(M,X)$ be a brick with
non-zero complexity. Let $P$ be a minimal skeleton of $(M,X)$, and let $Q$
be a filling of $P$. Then $S(Q)$ is connected and satisfies the following:

\begin{enumerate}

\item No pair of edges disconnects $S(Q)$.

\end{enumerate}
Suppose in addition either that every torus in $M$ is separating or
that $c(M,X)\le9$. Then:

\begin{enumerate}\addtocounter{enumi}{1}

\item If a quadruple of edges disconnects $S(Q)$, then one of the two
resulting components must be of one of the forms shown in
Fig.~\ref{2bridge}.

\end{enumerate} \end{teo}

An important tool of our search for bricks is the following
non-minimality criterion, proved in
Subsection~\ref{proof:shortloops}.
Let us say that a loop $\gamma$ in a skeleton $P$ of
$(M,X)\in\calX$ \emph{bounds an external disc} if there exists a closed
disc $D \subset M$ with $\partial D = \gamma$ and $D\cap P = \gamma$. A
loop is \emph{fake} if it is contained in the link of some point of $P$.

\begin{teo} \label{shortloops} Let $P$ be a standard skeleton of a
manifold with triods. If $P$ contains a non-fake loop which
bounds an external disc and has length at most $3$,
then $P$ is not minimal. \end{teo}

\paragraph{Computer search for bricks.} To find $\calB_{\le
9}\setminus\{B_0,\ldots,B_4\}$ we have first listed by computer the
4-valent graphs satisfying the conditions of
Theorem~\ref{restrictions:on:graph}. For each such graph $\Gamma$,
using~\cite{Be-Pe}, we have then determined the standard spines $Q$
of closed manifolds with $S(Q)\cong\Gamma$ and satisfying the
conditions of Remark~\ref{real:meaning:loops}. Then we have tested
the $Q\setminus\calR(\calL(Q))$'s for minimality using
Theorem~\ref{shortloops}. The result has been a very short list of
skeleta, but actually not all of them were minimal, and some pairs of
them were minimal skeleta of the same element of $\calX$. To
eliminate non-minimal and duplicate skeleta we have therefore used
certain moves on polyhedra which are known to transform a skeleton
$P$ into another skeleton $P'$ of the same $(M,X)$. Namely, we have
used the Matveev-Piergallini move and some disc-replacement moves
involving discs of length at most 4 (see \cite{Mat} for definitions).
The result has been a list of minimal skeleta of pairwise distinct
elements of $\calX^{\rm pr}$, but a few non-bricks were still
present. To get rid of them we have used very technical extra
criteria (such as Theorem~\ref{traces:4} below). The fact that the
list of 30 elements eventually obtained cannot be further reduced, so
all its elements are indeed bricks, follows from the (easy) fact that
no element of the list is obtained via a sharp-assembling from the
other ones.

\begin{rem} \emph{The bound $c(M,X)\le 9$ in
Theorem~\ref{restrictions:on:graph} is definitely not sharp, and we
actually conjecture the theorem to be true for any complexity. Moreover,
if an assembling of some bricks is a manifold in which each torus
is separating, then the same happens in the
individual bricks. Therefore, if one wants to
search for closed atoroidal manifolds only, the search for bricks can be
restricted to those in which each torus is separating, to which the whole
of Theorem~\ref{restrictions:on:graph} applies.} \end{rem}

We explain now how Theorem~\ref{restrictions:on:graph} helps saving space
in the search for bricks, by ruling out most graphs as possible $S(Q)$'s.
Namely, let $\calK$ be the set of all 4-valent graphs, let $\calK_{\rm
brick}\subset\calK$ consist of all $S(Q)$'s where $Q$ is a filling of some
minimal skeleton of some brick, and let $\calH\subset\calK$ consist of the
graphs satisfying both the constraints of
Theorem~\ref{restrictions:on:graph}. We know that $\calK_{\rm
brick}\subset\calH\subset\calK$ (at least in complexity $\le 9$, or for
bricks in which all tori are separating).
Table~\ref{skeleta:table} lists, up to 10 vertices, the
number of elements of each of these sets, showing that indeed $\#\calH$ is
a lot smaller than $\#\calK$, and not so far from $\#\calK_{\rm brick}$.
We still have not determined the bricks with 10 vertices.

\begin{table}
\begin{center}
\begin{tabular} {|l|c|c|c|c|c|c|c|c|c|c|}
\hline Vertices & $1$ & $2$ & $3$ & $4$ & $5$ & $6$
& $7$ & $8$ & $9$ & $10$ \\ \hline \hline $\calK$ & $1$ & $2$ & $4$ &
$10$ & $28$ & $97$ & $359$ & $1635$ & $8296$ & $48432$ \\ $\calH$ &
$1$ & $1$ & $1$ & $2$ & $4$ & $11$ & $27$ & $57$ & $205$ & $1008$ \\
$\calK_{\rm brick}$ & $1$ & $1$ & $1$ & $2$ & $3$ & $1$ & $4$ & $9$ &
$13$ & $?$
\\ \hline
\end{tabular}
\nota{Graphs which are singular sets of bricks.}
\label{skeleta:table}
\end{center}
\end{table}

\subsection{The non-minimality criterion} \label{proof:shortloops}
We prove here Theorem~\ref{shortloops}.

\begin{rem} \label{disc-rem} \emph{Let $(M,X)\in\calX$ be given together
with a standard skeleton $P$. A closed surface $F\subset P$ contains a
graph $H = F\cap S(P)$ with vertices having valency $3$ and $4$. If $F$ is
orientable, then we can choose a transverse orientation and give each edge
$e$ of $H$ a red or black color, depending on whether $P$ locally lies on
the positive or on the negative side of $F$ near $e$. A vertex with
valency $3$ is adjacent to edges with the same color, and a vertex with
valency $4$ is adjacent to two opposite red edges and two opposite black
ones.} \end{rem}

\medskip\noindent \emph{Proof of Theorem~\ref{shortloops}.}
Let $D$ be an external disc bounded by a loop
as in the statement. If we add
$D$ and remove a face in $\Sigma(P\cup D)$ we get another skeleton of $P$.
We prove now that there is a face in
$\Sigma(P\cup D)$ incident to more than $l(\partial D)$ distinct vertices.
This shows that $P$ is not minimal.

We consider the graph $H=S(P\cup D)\cap \Sigma(P\cup D)$, which contains
$\partial D$. By Proposition~\ref{loops-prop}, the surface $\Sigma(P\cup
D)$ is orientable; we can then choose a transverse orientation and color
the edges as explained in Remark~\ref{disc-rem}. Suppose by contradiction
that each face $f \subset \Sigma(P\cup D)$ meets at most $l(\partial D)$
vertices.

\begin{figure}
\begin{center}
\includegraphics{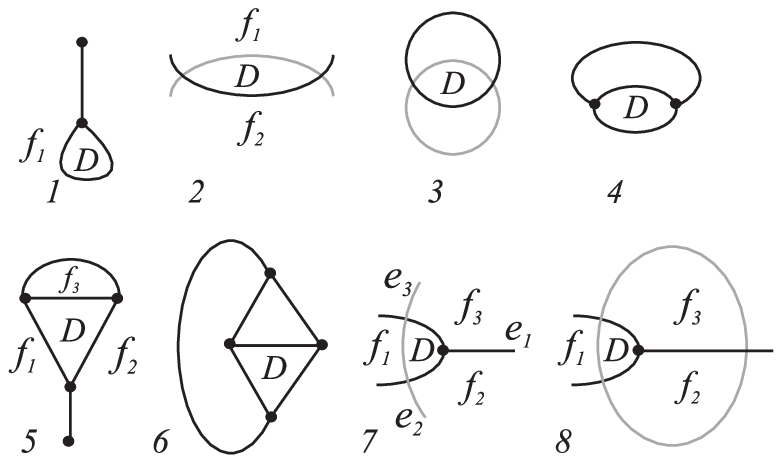}
\nota{Possible configurations
for $\Sigma(P\cup D)$ in the proof of Theorem~\ref{shortloops}.}
\label{sferette}
\end{center}
\end{figure}

A vertex in $\partial D$ has valency $4$ if and only if it is
adjacent to two distinct edges in $\partial D$ with distinct colors.
If $l(\partial D)=1$, then the only vertex contained in $\partial D$
would have valency $3$, as in Fig.~\ref{sferette}-($1$). So $f_1$
would meet at least $2$ distinct vertices.

If $l(\partial D)=2$, then the two vertices adjacent to $\partial D$ have
the same valency. Suppose they both have valency $4$, as in
Fig.~\ref{sferette}-($2$). Since each $f_i$ meets at most $2$ vertices,
then $H$ is as in Fig.~\ref{sferette}-($3$), but $M\setminus(P\cup D)$
would have 3 components. Suppose both vertices of $D$ have valency $3$:
then $H$ is as in Fig.~\ref{sferette}-($4$), and $\partial D$ is fake.
Both cases are excluded.

If $l(\partial D)=3$, either all vertices met by $\partial D$ have valency
$3$, or two of them have valency $4$. Suppose the first case holds. If a
face $f_i$ meets $2$ distinct vertices only, then the other two faces
adjacent to $D$ coincide, as in Fig.~\ref{sferette}-($5$), and meet more
than $3$ vertices. So each $f_i$ meets $3$ distinct vertices, and $H$ is
the 1-skeleton of the tetrahedron $\Sigma(P\cup D)$ as in Fig.~\ref{sferette}-(6);
hence $\partial D$ is fake, which is absurd.

Finally, suppose two vertices have valency $4$ and one has valency $3$ as
in Fig.~\ref{sferette}-(7); since $f_2$ is incident to at most $3$
distinct vertices, then the distinct edges $e_1, e_2$ have one common
endpoint; for the same reason the edges $e_1, e_3$ have one common
endpoint. Then $H$ is as in Fig.~\ref{sferette}-(8); but this is absurd
since $M\setminus(P\cup D)$ would have at least $3$ components.


\vspace{-.85cm}
\begin{flushright} $\square$ \end{flushright}

\section{Minimal skeleta of bricks}\label{proofs:section}

In this section we will prove
Theorems~\ref{minimal:are:standard},~\ref{disjoint:length:1:loops},
and~\ref{restrictions:on:graph}. This requires the introduction of several
ideas not mentioned yet. The crucial point of our work will be the
analysis of the intersection between a minimal skeleton and a closed
orientable surface. We warn the reader that the proofs of
Theorems~\ref{traces:2:difficult} and~\ref{traces:4} given below are long
and not very much illuminating by themselves, so they can be skipped at
first. We will only consider in this section bricks having {\em positive}
complexity, without further notice. So all minimal skeleta will be standard
by Theorem~\ref{standard:teo}.

\subsection{Traces}

\begin{figure}
\begin{center}
\includegraphics{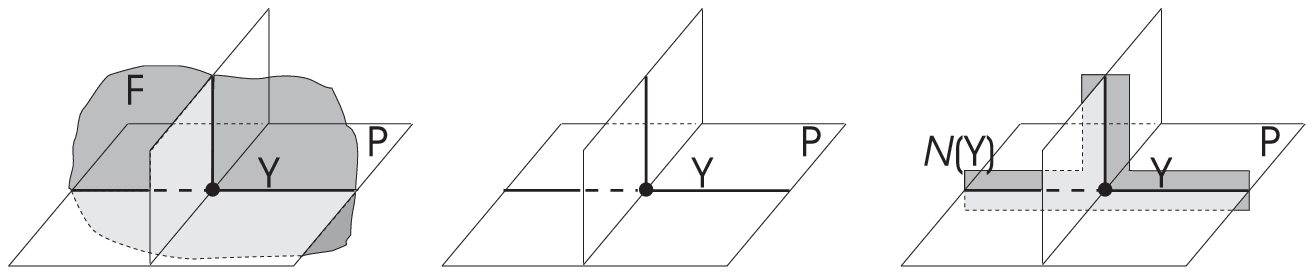}
\nota{Surface with boundary $\calN(Y)$
constructed from the trivalent graph $Y$.}
\label{st-fig} \end{center}
\end{figure}

Let $(M,X)$ be a manifold with triods and let $P$ be a standard skeleton
of $(M,X)$. A closed surface $F \subset \interior M$ is said to be
\emph{simply transverse} to $P$ if:

\begin{enumerate}

\item $F$ is transverse to $P$;

\item The intersection of $F$ with $M\setminus P$ consists of a finite
number of discs.

\end{enumerate} The first condition implies that $Y=P\cap F$ is a finite
trivalent graph disjoint from $V(P)$, whose vertices lie precisely at the
intersection of $Y$ with the edges of $P$ and appear as in
Fig.~\ref{st-fig}-left. Such a graph is called the \emph{trace} of $F$.

\begin{rem} \emph{Let a trivalent graph $Y\subset P \setminus V(P)$ be
given, in such a way that $Y \cap S(P)$ consists of all the vertices of
$Y$, each appearing as in Fig.~\ref{st-fig}-centre. We show that $Y$ is
the trace of an essentially unique simply transverse surface $F \subset
M$. First, we can uniquely find a surface $\calN(Y)$ with boundary,
transverse to $P$, which collapses to $Y$ (see $\calN(Y)$ near an edge of
$P$ in Fig.~\ref{st-fig}-right). The boundary of $\calN(Y)$ consists of a
finite number of circles that lie on the boundary of a sub-ball $B'$ of
$B$. Now we can uniquely glue disjoint discs properly embedded in $B'$ to
these circles, thus getting the desired closed surface $F$.} \end{rem}

\subsection{Traces with 2 vertices} \label{traces:2}

\begin{lemma} \label{properties-cor2} Let $(M,X)$ be a brick and let $P$
be a minimal skeleton of $(M,X)$. Let $Y\subset P$ be the trace of an
orientable surface $F\subset M$. Then each edge of $Y$ has distinct
endpoints. \end{lemma}

\begin{proof} Suppose $s$ is an edge of $Y$ with common endpoints; since
$F$ is orientable, the regular neighbourhood of $s$ in $F$ is an
annulus, so there is a component $D_0 \subset F\setminus Y$ with
$\partial D_0=s$. Then $\partial D_0$ is a length-1 loop; this is
impossible by Theorem~\ref{shortloops}, since length-1 loops are
never fake. \end{proof}

Let $P$ be a standard skeleton of some $(M,X)\in\calX$ and let $\theta_i
\in X$ be the triod contained in $T_i\subset\partial M$. Pushing $\theta_i$
a bit inside $\interior P$ we get the trace $Y$ of a torus parallel to
$T_i$. Therefore we say that such a $Y$ is \emph{parallel to the boundary}
(of $P$).

\begin{prop}\label{traces:2:easy} Let $(M,X)$ be a brick, equipped with a
minimal skeleton $P$. Let $Y$ be the trace of an
orientable surface $F$. If $Y$ has two vertices
then it is a triod and one of the following
occurs:

\begin{enumerate}

\item $F$ is a non-separating torus;

\item $Y$ is parallel to the boundary;

\item $F$ is a sphere and $Y$ is the link of a point in $S(P)\setminus
V(P)$.

\end{enumerate} \end{prop}

\begin{proof} First, $Y$ is a triod by Lemma~\ref{properties-cor2}. There
are two possibilities for the regular neighbourhood $\calN(Y)$ of $Y$
in $F$, which are shown in Fig.~\ref{triods} and lead to a sphere and
a torus respectively. In the first case $F\setminus Y$ contains three
external discs $D_i$ with $e(D_i)=2$. By Theorem~\ref{shortloops} all
the loops $\partial D_i$ are fake, so $Y = {\rm lk}(p)$ for some $p
\in S(P)\setminus V(P)$.
\begin{figure}
\begin{center}
\includegraphics{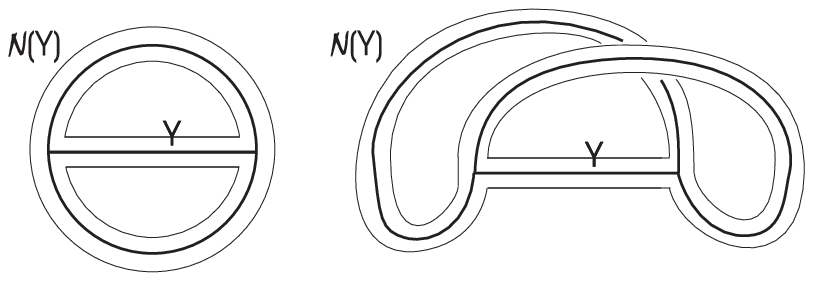}
\nota{The two possibilities for an orientable $\calN(Y)$ when $Y$ is a triod.} \label{triods}
\end{center} \end{figure}
In the second case, let $F$ be separating, and let $N_1$ and $N_2$ be the
manifolds into which $F$ separates $M$. Set $P_i = N_i\cap P$ for $i=1,2$.
Then $(M,X)$ is obtained by assembling the manifolds with triods
$(N_1,X_1)$ and $(N_2,X_2)$, where $X_i = \cup\partial P_i$ for $i=1,2$.
Moreover $P_i$ is a skeleton of $(N_i,X_i)$, which implies that this
assembling is sharp unless it is trivial. Since $(M,X)$ is a brick,
the assembling is trivial. Now, $P_1$ and $P_2$ are standard, so
$(M_1,X_1)$ and $(M_2,X_2)$ are prime by Theorem~\ref{standard2:teo}.
Therefore, the assembling must be of the first trivial type, namely
$(N_1,X_1)$ must be  $B_0$ up to permutation.
Hence $P_1$ is the unique minimal skeleton of $B_0$,
homeomorphic to $\theta\times [0,1]$. It follows that $Y$ is parallel to
the boundary in $P$. \end{proof}

\begin{cor} \label{no:embedded:faces} Let $P$ be a minimal skeleton of a
brick. Then there is no embedded face  in $P$ incident to $3$ or fewer
vertices. Moreover, for every edge $e$ of $P$, the three faces of $P$ incident
to $e$ are distinct from each other. \end{cor}

\begin{proof} Let $f$ be an embedded face with $3$ or fewer vertices. A
loop in $P$ very close to $\partial f$ and disjoint from
$\overline{f}$ bounds a disc $D$ parallel to $f$. Moreover
$l(\partial D)\le 3$ and $\partial D$ is not fake since
$M\setminus(P\cup\partial M)$ has only one component.

Let $f\subset P$ be a face incident at least twice to an edge $e$ of
$P$. It follows that there is a length-1 loop $\gamma\subset P$
intersecting $e$ once. Length-1 loops are never fake, so, by
Theorem~\ref{shortloops}, $\gamma$ does not bound a disc. Therefore
its regular neighbourhood $\calR(\gamma)$ is a M\"obius strip with
one tongue, and $\partial \calR(\gamma)$ is a trace with two vertices
of the disconnecting torus in $M$ which bounds the regular
neighbourhood of $\gamma$ in $M$. Proposition~\ref{traces:2:easy}
implies that $\partial\calR(\gamma)$ is boundary-parallel, so $P$ has
no vertices.
\end{proof}

\paragraph{Co-disconnecting surfaces.} Let $P$ be a standard skeleton of
$(M,X)\in\calX$. Let $Y\subset P$ be the trace of a simply transverse
orientable surface $F$. Every component $D$ of $F\setminus Y$ is an open
disc; its boundary is the union of two parts $\partial_1 D$ and
$\partial_2 D$, where $\partial_i D$ is the closure of the union of all
edges of $Y$ adjacent $i$ times to $D$.
If we add $D$ to $P$ we do not get a simple polyhedron, unless $\partial_2
D =\emptyset$. It is nevertheless easy to see that
Proposition~\ref{loops-prop} holds for $P\cup D$ too, namely:

\begin{prop} \label{loops-prop2} Let $B'$ and $B''$ be the balls given by
$M\setminus (P\cup D)$. The faces of $P\cup D$ that separate $B'$ from
$B''$ form a closed orientable surface $\Sigma(P\cup D) \subset P\cup D
\subset M$ which cuts $M$ into two components. \end{prop}

\begin{proof} The proof of Proposition~\ref{loops-prop} works away from
$\partial_2 D$. We only need to show that $\Sigma(P\cup D)$ is a surface
near $\partial_2 D$: let $f'$ and $f''$ be the faces other than $D$
incident to an edge $e\subset\partial_2 D$. Since $F$ is orientable, $f'$
is adjacent to $B'$ on both sides and $f''$ is adjacent to $B''$ on both
sides (or the converse). Therefore $f'$ and $f''$ are disjoint from
$\Sigma(P\cup D)$, and $\Sigma(P\cup D)$ is a closed surface. \end{proof}

In the above setting we define $\Sigma_D\subset P$ as $\Sigma(P\cup
D)\setminus D$, and call it the \emph{co-disconnecting surface} of
$D$. By Proposition~\ref{loops-prop2}, $\Sigma_D$ is a compact
surface with boundary $\partial_1 D$. For a subpolyhedron $K\subset
P$ we will denote from now on by $\calR(K)$ and $\calR_M(K)$ the
regular neighbourhoods of $K$ in $P$ and in $M$ respectively.

\begin{teo}\label{traces:2:difficult} Under the assumptions of
Proposition~\ref{traces:2:easy}, assume that $c(M,X)\le 9$. Then $Y$
cannot be the trace of a non-separating torus. \end{teo}

\begin{proof} By contradiction let $Y=T\cap P$ with $T$ non-separating,
and put $D=T\setminus Y$. The co-disconnecting surface $\Sigma_D\subset P$
is by Proposition~\ref{loops-prop2} a closed orientable surface, which is
non-empty since $\Sigma_D\cup T$ disconnects $M$, whereas $T$ does not. We
assume that $\#(V(P)\cap\Sigma_D)$ is minimal among all mimimal skeleta of
$(M,X)$ for which there exists a non-separating torus whose trace is a
triod.

We focus now on a component $\Sigma$ of $\Sigma_D$. Choosing a transverse
orientation for $\Sigma$ as in Remark~\ref{disc-rem}, we can trace on
$\Sigma$ two trivalent graphs $Y_+$ and $Y_-$ which intersect
transversely. These graphs represent the way the rest of $P$ glues to
$\Sigma$, and the sign $+$ or $-$ depends on whether $P$ locally lies on
the positive or on the negative side of $\Sigma$. We show now several
properties of the triple $(\Sigma,Y_+,Y_-)$ which do not require the bound
9 on complexity. Only later we will use this bound.

\begin{enumerate}

\item\label{mainteo:planar:compo:point} \emph{$\Sigma\setminus Y_\pm$
consists of planar surfaces}. Given a point $p$ of
$\Sigma\setminus(Y_+\cup Y_-)$ there are two points $p_+,p_-$ of $\partial
\calR_M(P)$ closest to $p$, with $p_+$ on the positive side of $\Sigma$
and $p_-$ on the negative side. It is not hard to show that the map
$p\mapsto p_+$ extends to a homeomorphism of $\Sigma\setminus Y_+$ onto an
open subset of $\partial \calR_M(P)\cong S^2$, and similarly for $Y_-$.

\item\label{mainteo:bounding:circles:point} \emph{The components of
$\partial\calR(Y_\pm)\cap\Sigma$ bound discs in $M$.} This follows from the same
argument just explained.

\item\label{mainteo:disc:compo:point} \emph{$\Sigma\setminus(Y_+\cup Y_-)$
consists of discs.} This is because $\Sigma\subset P$, $Y_+\cup
Y_-=\Sigma\cap S(P)$, and $P$ is standard.

\item\label{mainteo:essential:circles:point} \emph{If a component of
$\Sigma\setminus \calR(Y_\pm)$ is not a disc then its boundary loops are
essential in $\Sigma$.} We refer to $Y_+$. If one of them is not, it is
very easy to see that there is a disc $\Delta$ in $\Sigma$ such that
$Y_+\cap\partial\Delta=\emptyset$ but $Y_+\cap\Delta\neq\emptyset$, so in
particular $Y_+\cap\Delta$ contains vertices of $P$. The move suggested in
Fig.~\ref{liftdisc} then contradicts minimality of $P$.

\begin{figure}
\begin{center}
\includegraphics{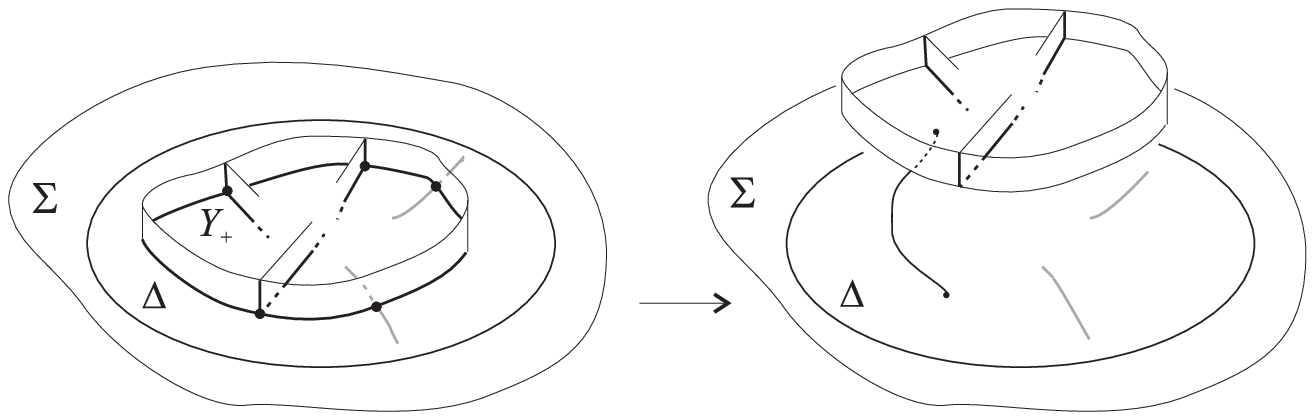}
\nota{A move which reduces complexity.}
\label{liftdisc}
\end{center}
\end{figure}

\item\label{mainteo:non-planar:R:point} \emph{Not all the components of
$\calR(Y_\pm)\cap\Sigma$ are planar.} Again we refer to $Y_+$. By contradiction,
from points~\ref{mainteo:planar:compo:point}
and~\ref{mainteo:bounding:circles:point} and the irreducibility of $M$, we
would readily get that $\Sigma$ bounds a handlebody, but $\Sigma$ is
non-separating.

\item\label{mainteo:1:intersection:point} \emph{Every component of $Y_+$
intersects $Y_-$, and conversely.} Otherwise, since $\Sigma$ is connected,
there would exist a component of $\Sigma\setminus(Y_+\cup Y_-)$ with
disconnected boundary, contradicting point~\ref{mainteo:disc:compo:point}.

\item\label{mainteo:2:intersections:point} \emph{$Y_+\cap Y_-$ contains at
least two points.} Assume there is only one point $v$ (a crossing between
$Y_+$ and $Y_-$). If a face $f$ of $\Sigma$ is incident to $v$, then it
must be multiply incident, because faces contain an even number of
quadrivalent vertices (with multiplicity). If two instances of $f$ are
adjacent to each other at $v$, we find in the closure of $f$ a length-1
loop bounding an external disc, which contradicts minimality. If two
instances of $f$ are opposite at $v$, then for the same reason there is
another face $g$ doubly incident to $v$, and $g\neq f$. Now in the closure
of $f\cup g$ we can easily find a length-2 loop bounding an external disc
which meets edges opposite at $v$. By minimality the loop must be fake, so
these edges must actually be the same. Orientability of $\Sigma$ then
implies that $f=g$: a contradiction.

\item\label{mainteo:many:intersections:point}
\emph{If a component of $Y_+$ is a circle then it intersects $Y_-$ in
at least 4 points, and conversely.} This readily follows from
Corollary~\ref{no:embedded:faces} and
minimality, because this circle is precisely the boundary of a face of $P$.

\item\label{mainteo:forbidden:squares:point} \emph{No squares as in
Fig.~\ref{badsquare}-left occur in $(\Sigma,Y_+,Y_-)$.} If one such square
exists, we can correspondingly apply to $P$ one move as in
Fig.~\ref{badsquare}-right. The result is a new
minimal skeleton $P'$ on which $T$ still has trace $Y$, but
$\#(V(P')\cap\Sigma'_D)<\#(V(P)\cap\Sigma_D)$. A contradiction.
\begin{figure}
\begin{center}
\includegraphics{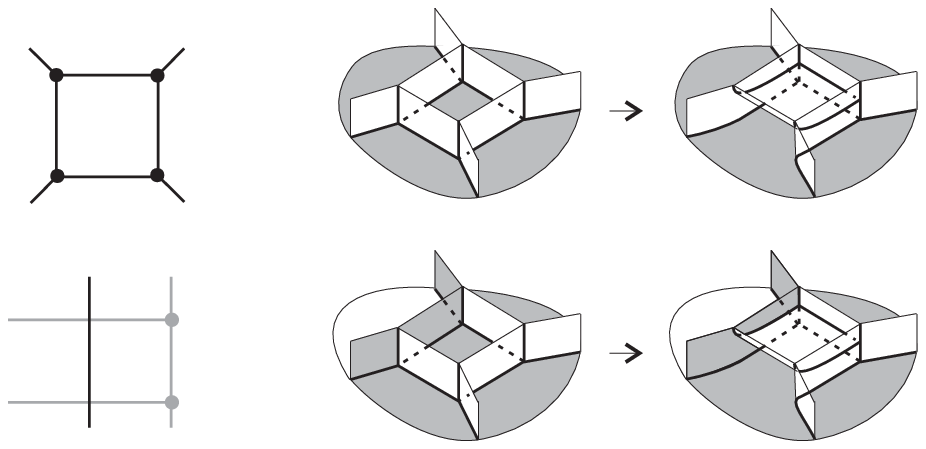}
\nota{Forbidden squares (with black $Y_+$ and gray $Y_-$),
and moves (with shadowed $\Sigma$).} \label{badsquare}
\end{center}
\end{figure}

\end{enumerate}

We show now how to conclude, using the fact that $\#V(P)\le 9$. It follows
from point~\ref{mainteo:non-planar:R:point} that both $Y_+$ and $Y_-$ have
vertices. Being trivalent, they have an even number of them, and the total
is at most $9-2=7$ by point~\ref{mainteo:2:intersections:point}. So up to
permutation we can assume that $Y_+$ has 2 vertices. In particular
$\calR(Y_+)\cap\Sigma$ has only one non-planar component, which is homeomorphic to a
punctured torus (with a component of $Y_+$ sitting as a triod in this
torus). From point~\ref{mainteo:many:intersections:point} we deduce that
$Y_+$ can have at most one circular component, and it is now easy to
deduce from point~\ref{mainteo:planar:compo:point} that $\Sigma$ indeed is
a torus. Point~\ref{mainteo:essential:circles:point} then implies that
$Y_+$ consists of the triod only. In the rest of our proof we will always
depict $\Sigma$ cut open along $Y_+$. So $\Sigma$ appears as a hexagon,
and we denote by $\Delta$ its interior.

To conclude the proof we will first show that $Y_-$ also has 2 vertices,
and then that it appears in one of the two shapes shown in
Fig.~\ref{goody}. This indeed yields a contradiction to the fact that
$(M,X)$ is a brick, since $Y_-\cap Y_+$ consists of two points, so cutting
$P$ along $\Sigma$ we see that $(M,X)$ can be obtained via a
sharp-assembling.
\begin{figure}
\begin{center}
\includegraphics{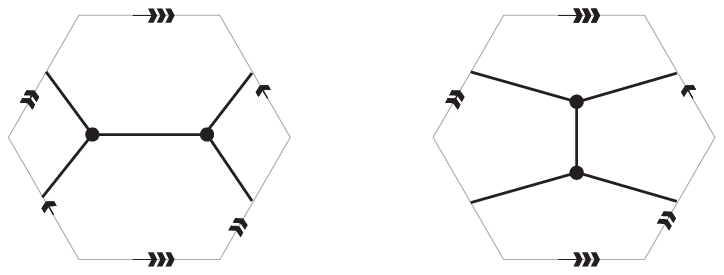}
\nota{Configurations corresponding to a self-assembling.}
\label{goody}
\end{center}
\end{figure}

So, let $Y_-$ have 4 vertices. We claim that all the components of
$Y_-\cap\Delta$ are trees. If one of them is not then there is a face
of $P$ inside $\Delta$ and bounded by $Y_-$. Then either this face
has $\le 3$ vertices, which contradicts
Corollary~\ref{no:embedded:faces}, or it is a square of the first
forbidden type. Our claim is proved. Now note that if $Y_-\cap\Delta$
has $\nu$ components then it has $4+2\nu$ free endpoints, which give
$2+\nu$ vertices in $P$. Since $Y_+$ has 2 vertices and $Y_-$ has 4,
we deduce that $\nu=1$ and that $Q=P\setminus\calR_P(\Sigma)$ has no
vertices. Moreover $Q$ is connected and standard, and $\partial
Q\cong Y_+\sqcup Y_-$. It is not hard to show that with these
constraints the only possibility for $Q$ is as shown in
Fig.~\ref{no4verts},
\begin{figure}
\begin{center}
\includegraphics{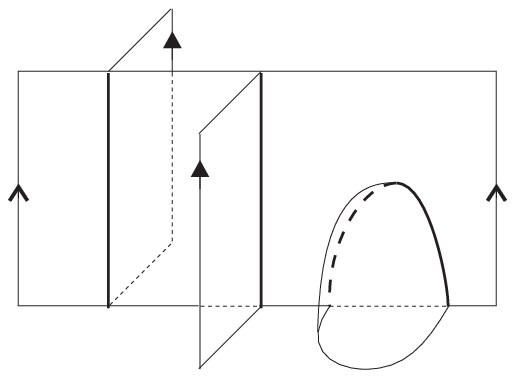}
\nota{A polyhedron without vertices.}
\label{no4verts}
\end{center}
\end{figure}
so $\partial\calR_M(Q)$ has two components. In addition, also
$\Sigma\setminus Y_-$ consists of discs (as $\Sigma\setminus Y_+$), and we
get a contradiction because $\partial\calR_M(Q)$ should then be a sphere
with some holes.

We can now assume that $Y_-$ has two vertices, and show that it appears as
in Fig.~\ref{goody}. Knowing already that $\Sigma\setminus Y_-$ is a disc,
it is enough to show that $Y_- \cap\Delta$ is connected. Suppose by
contradiction that $Y_- \cap\Delta$ is disconnected. Then there exists an
arc $\alpha$ properly embedded in $\Delta$ which separates two components
of $Y_- \cap \Delta$. Let us consider the endpoints of $\alpha$. By
minimality of $P$, they cannot belong to the same edge of $\Delta$, nor to
two adjacent ones, otherwise we could make $Y_-$ slide on $\Sigma$ and
reduce the number of vertices, as in Fig.~\ref{connect1}.
\begin{figure}
\begin{center}
\includegraphics{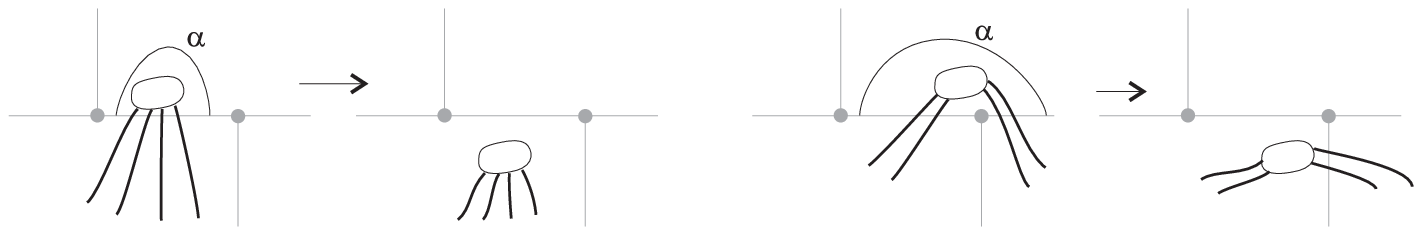}
\nota{Moves reducing complexity.}
\label{connect1}
\end{center}
\end{figure}
The ends of $\alpha$ also cannot belong to two edges adjacent to one
and the same edge, as in Fig.~\ref{connect2}-left.
\begin{figure}
\begin{center}
\includegraphics{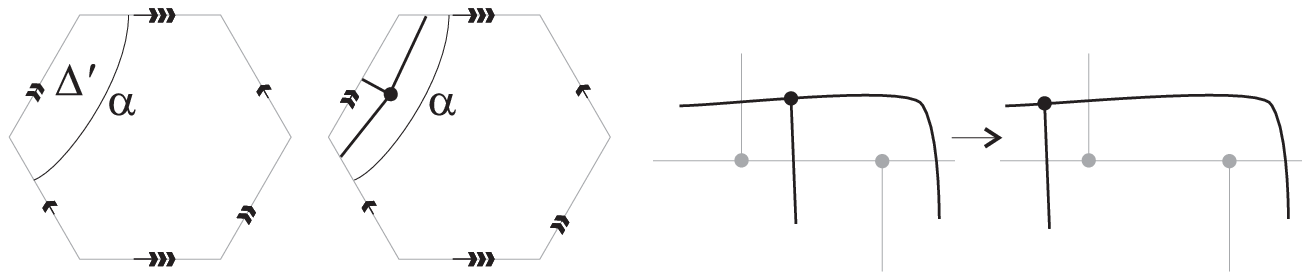}
\nota{More moves reducing complexity.}
\label{connect2}
\end{center}
\end{figure}
To see this, consider how many vertices of $Y_-$ can lie in $\Delta'$. If
there are no vertices at all, then either a face of $P$ contained in
$\Delta'$ has less than 4 vertices or there is a square of the second
forbidden type.  If $Y_-$ has both vertices in $\Delta'$, then again
$\Delta'$ contains either a small face or a forbidden square. These cases
are excluded, so there is one vertex of $Y_-$ in $\Delta'$, and the only
possible case is shown in Fig.~\ref{connect2}-center. Now we let $Y_+$
slide over $\Sigma$ as shown in Fig.~\ref{connect2}-right. The result is a
new minimal skeleton $P'$ on which $T$ still has trace $Y$, but
$\Sigma'_D$ now contains one of the forbidden squares of
Fig.~\ref{badsquare}, which contradicts minimality of
$\#(V(P)\cap\Sigma_D)$.

We are left to show that the endpoints of $\alpha$ also cannot belong to
opposite edges of $\Delta$ (Fig.~\ref{connect3}-left).
\begin{figure}
\begin{center}
\includegraphics{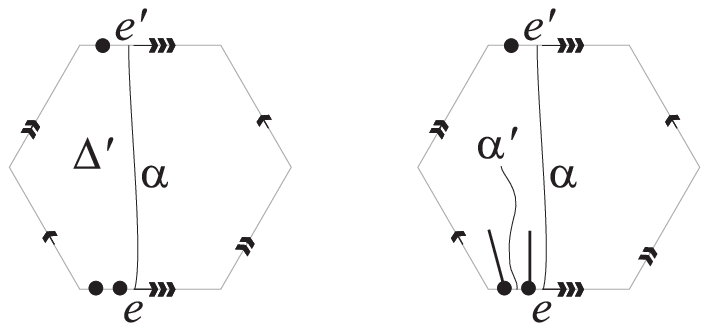}
\nota{Conclusion of Step 3.}
\label{connect3}
\end{center}
\end{figure}
Denote by $\nu$ and $\nu'$ the number of ends of $Y_-\cap\Delta$ on $e\cap
\Delta'$ and on $e'\cap\Delta'$ respectively. If $\nu=\nu'$ then $\alpha$
can be isotoped so to give rise to a length-1 loop in $P$ bounding an
external disc: a contradiction. If $\nu=0$ or $\nu'=0$ then we can replace
$\alpha$ by a curve disjoint from $Y_-$ and having ends on edges of
$\Delta$ which are not opposite, so we get back to a case already ruled
out. So up to permutation we can assume that $\nu\geq 2$. Now the face of
$P$ containing the portion of arc $\alpha'$ shown in
Fig.~\ref{connect3}-right must meet another edge of $Y_+=\partial\Delta$,
otherwise it is either small or forbidden (recall that $Y_-$ has 2
vertices only). So $\alpha'$ extends to a properly embedded arc disjoint
from $Y_-$. Either $\alpha'$ belongs to a case already ruled out, or the
corresponding $\nu+\nu'$ is smaller, and a contradiction is reached
anyway. This eventually shows that $Y_-$ is connected, and the proof is
complete. \end{proof}

\subsection{Moves on traces} \label{J}

\begin{figure}
\begin{center}
\includegraphics{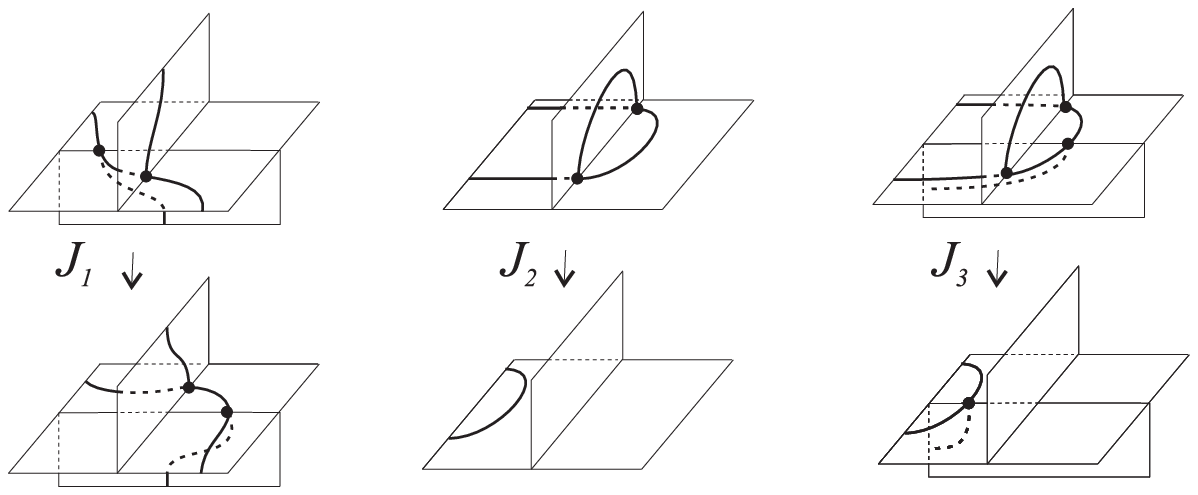}
\nota{The moves $J_1, J_2, J_3$.} \label{j}
\end{center} \end{figure}

The key step to check the properties of bricks will be
Theorem~\ref{traces:4} stated below. We introduce here more new notions
which will be used to prove it.

Let $P$ be a standard skeleton of a manifold with triods $(M,X)$. Given
the trace $Y$ of a surface $F$, there are some obvious moves that
transform $Y$ into another trace $Y'$ of a surface $F'$ isotopic to $F$.
Three such moves, denoted by $J_1$, $J_2$ and $J_3$ and collectively
called $J$-moves, are shown in Fig.~\ref{j}. Since we will be concerned
with traces of (transversely) orientable surfaces only, we can ask a
$J$-move to transform a trace $Y$ into a trace $Y'$ disjoint from $Y$. Let
$[Y,Y']$ be the sub-polyhedron which lies between $Y$ and $Y'$. A sequence
of moves $Y_1 \to \ldots \to Y_n$ is called a \emph{flow} if each move
$Y_i\to Y_{i+1}$ is a $J$-move and $[Y_{i-1},Y_i] \cap [Y_i, Y_{i+1}] =
Y_i$ for all $i$, namely, if the moves are performed towards the same
normal direction to $Y_i$ for all $i$.

\begin{rem}\label{J1:along} \emph{A move $J_1$ is determined by an edge
$s$ of $Y$ and a vertex $v$ of $P$ such that $s\subset{\rm lk}(v)$, or
equivalently by the cone $vs$ from $v$ based on $s$ (a triangle). We will
sometimes say that the move is performed \emph{along} the triangle.}
\end{rem}

\begin{rem} \emph{If a move $J_1$ transforms a trace $Y$ of $F$ into a
trace $Y'$ of $F'$ then there is a natural bijection between the
components of $F\setminus Y$ and those of $F'\setminus Y'$.} \end{rem}

Let $Y$ be the trace of a surface $F$. Given a component $D$ of
$F\setminus Y$, we denote by $e(D)$ the number of edges of $Y$ adjacent to
$D$, counted with multiplicity (\emph{i.e.}~an edge of $Y$ is counted
twice if it has $D$ on both sides).

\begin{lemma} \label{J-lemma:2} Let $P$ be a minimal skeleton of a brick.
Let $Y\subset P$ be the trace of an orientable surface $F$, and let $D$ be a
component of $F\setminus Y$. Consider a move $J_1$ determined by an edge
$s\subset\partial D$ of $Y$ and a vertex $v$ of $P$, call $Y'$ the
resulting trace and $D'$ the disc corresponding to $D$. Then $e(D') <
e(D)$ if $e(D)<6$ and $e(D') \le e(D)$ if $e(D)=6$. \end{lemma}

\begin{proof} The trace $Y'$ is obtained from $Y$ as shown in
Fig.~\ref{bound6}; it follows from the figure that if $e(D')>e(D)$
\begin{figure}
\begin{center}
\includegraphics{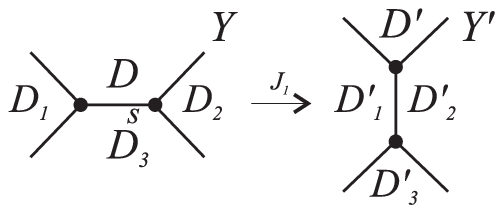}
\nota{A move $J_1$ at the level of traces.}
\label{bound6} \end{center} \end{figure}
then $D_1=D_2=D\neq D_3$ and if $e(D')=e(D)$ then $D=D_1$ or $D=D_2$. By
Lemma~\ref{properties-cor2} the edges of $Y$ have distinct ends. Using
this fact one easily sees that $e(D)>6$ if $D_1=D_2=D\neq D_3$ and
$e(D)\ge 6$ if $D=D_1$ or $D=D_2$, and the conclusion follows. \end{proof}

\paragraph{Good discs.} Let $Y$ be the trace of a surface $F$. We say that a
disc $D\subset F \setminus Y$ is \emph{good} if all discs in $F\setminus
Y$ other than $D$ are contained in the same component of $M\setminus
(P\cup D)$.

\begin{rem} \label{J-lemma:1} \emph{If $F$ has 2 discs then these discs
are good.} \end{rem}

\begin{rem}\label{ultimo:ritocco?} \emph{If $F$ is orientable, then
$\calR_P(Y)\cong Y\times [-1,1]$. Recall that
$\partial\Sigma_D=\partial_1D$. Now it is not hard to show that if $D$ is
good then the identification $\calR_P(Y)\cong Y\times [-1,1]$ can be
chosen so that $\Sigma_D \cap \calR_P(Y)\cong\partial_1 D \times [0,1]$,
and the converse holds if $F$ is connected. In other words, when $F$ is
orientable and connected, we have that $D$ is good if and only if
$\Sigma_D$ lies on a definite side of $Y$ in $P$.} \end{rem}

\begin{lemma} \label{J-lemma:3} Under the assumptions of
Lemma~\ref{J-lemma:2}, suppose that $D$ is good and that $\Sigma_D$ and the
triangle $vs$ lie on the same side of $Y$ in $P$. Then $D'$ is good and
$\Sigma_{D'} = \Sigma_D\setminus [Y,Y']$. \end{lemma}

\begin{proof} The condition that $\Sigma_D$ and $sv$ lie on the same of
side of $Y$ means that $Y$, during its transformation into $Y'$, is pushed
towards $\Sigma_D$, and the conclusion is obvious. \end{proof}

\subsection{Traces with 4 vertices}

We prove here the key result needed to
establish the properties of bricks.

\begin{rem} \label{pictures} {\em If a polyhedron $P$ is super-standard
(with boundary), then it can be uniquely reconstructed from the
regular neighbourhood $\calR(S(P))$ of $S(P)$ in $P$, by gluing discs
to each circle in $\partial\calR(S(P))$, because the rest of
$\partial\calR(S(P))$ can be identified to $\partial P$. Therefore
here and in the sequel we will describe such $P$'s by drawing
$\partial\calR(S(P))$ in $\matR^3$. Three-dimensional pictures will
be needed when $P$ is only standard.}
\end{rem}

\begin{teo} \label{traces:4} Let $P$ be a minimal skeleton of a brick
$(M,X)$, and let $Y\subset P$ be the trace of an orientable connected
surface $F$ with $4$ vertices. If $F$ is separating, then $Y$ is a
boundary component of a polyhedron of one of the following types:

\begin{enumerate}

\item\label{vert:nhb:item} ${\cal R}_P(v)$ for some $v\in V(P)$ (type
1.1), or $\calR_P(\lambda)$ for an arc $\lambda$ properly embedded in a
face of $P$ (type 1.2);

\item\label{loop:item} $\calR_P(\gamma)$ for a length-2 loop $\gamma$,
which is fake if it bounds an external disc;

\item\label{standard:item} One of the 5 polyhedra shown in
Fig.~\ref{type3}, whose boundary has two components: $Y$ and a triod
$\theta_i\subset\partial P$;

\begin{figure} \begin{center}
\includegraphics{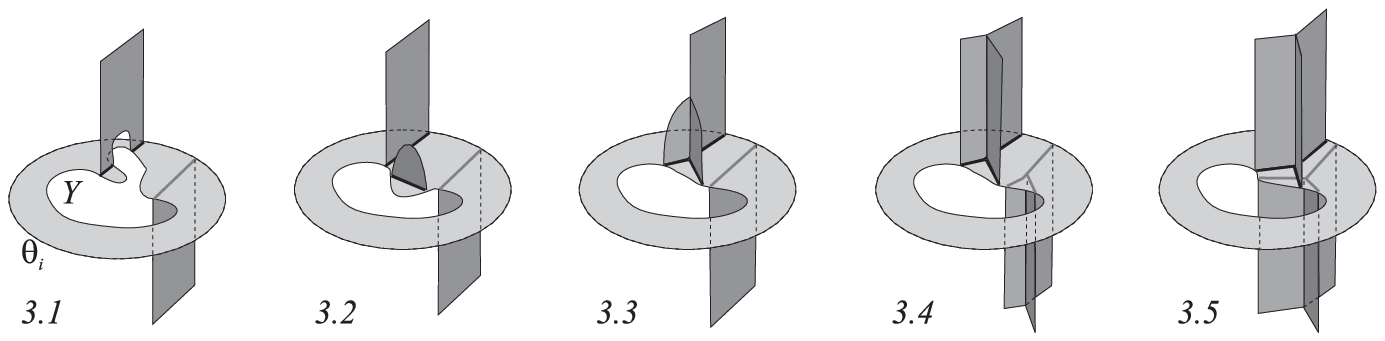}
\nota{Polyhedra of type 3: segments (in the first 3 pictures) and $Y$'s on the top are
identified to the corresponding ones on the bottom.} \label{type3} \end{center}
\end{figure}

\item\label{chain:item} A polyhedron as in Fig.~\ref{type4}, with $1$
(left) or more (right) vertices;

\begin{figure} \begin{center}
\includegraphics{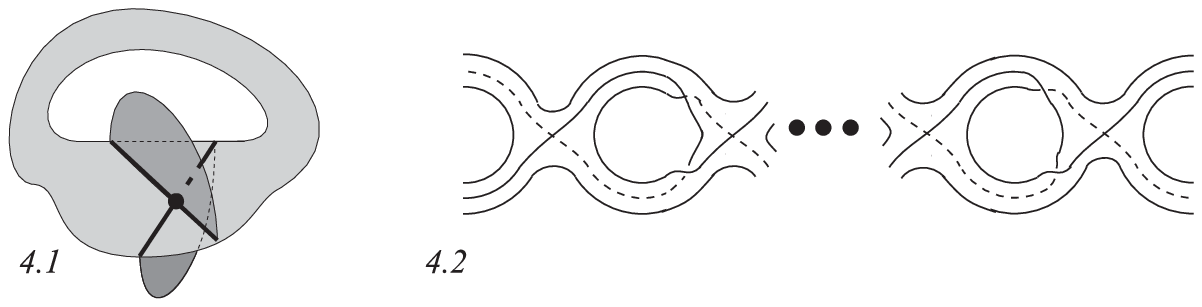}
\nota{Polyhedra of type 4.} \label{type4} \end{center} \end{figure}
\end{enumerate}
\begin{figure} \begin{center}
\includegraphics{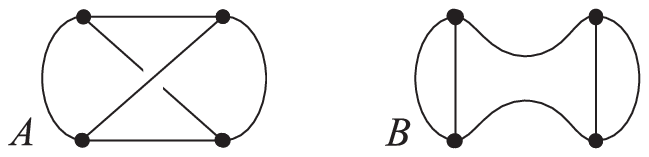}
\nota{Types $A$ and $B$ for $Y$.}
\label{abgraphs} \end{center}
\end{figure}

If $F$ is not separating, then it is a torus and there is a minimal
skeleton $P'$ of $(M,X)$ on which $F$ has a triod as a trace.

Moreover, only two types $A$ and $B$ of $Y$ are possible, as shown in
Fig.~\ref{abgraphs}. The polyhedra of types 1.1, 3.3, and 3.4 have
boundaries of type $A$, those of types 1.2, 3.1, 3.2, 3.5, and 4 have
boundaries of type $B$; a polyhedron of type \ref{loop:item} has boundary
of type $A$ if it is a M\"obius strip with two tongues, of type $B$
otherwise. \end{teo}

\begin{proof} It is enough to show that one of the following
must hold:
\begin{itemize}
\item[(I)] $F$ is a non-separating torus, and
$F$ has a triod as a trace on some $P'$;
\item [(II)] $Y$ bounds one of the polyhedra of type 1-4.
\end{itemize}
So we assume (I) does not hold and show (II). Our argument is long
and organized in many steps. We first describe the overall scheme
stating without proof 5 assertions. Later we will provide full
proofs. Let $D\subset F\setminus Y$ be a component having lowest
$e(D)$.

\smallskip\noindent{\bf Fact 1.} \emph{If $e(D) \in \{2,3\}$ then $Y$
bounds a polyhedron of type \ref{vert:nhb:item}, \ref{loop:item},
\ref{standard:item}.1, \ref{standard:item}.2, or \ref{standard:item}.3}
\smallskip

Suppose then that $e(D)\ge 4$. Since $Y$ is trivalent it has 6 vertices,
so $\chi(F)=d-2$, where $d$ is the number of components of $F\setminus Y$.
Each component is incident to at least 4 vertices, so $2\cdot 6\ge 4\cdot
d$, whence $d \le 3$. It easily follows that $F$ is a torus and $d=2$.
Then $F \setminus Y$ consists of two discs $D=D_1$ and $D_2$, both good by
Remark~\ref{J-lemma:1}. Recalling from Lemma~\ref{properties-cor2}
that all edges of $Y$ have distinct endpoints one easily sees that only
the types $A$ and $B$ for $Y$ are possible. The restriction that $e(D)\geq
4$ then implies that up to homeomorphism there is only one possible
configuration $(F,Y_A)$ and only one $(F,Y_B)$, as shown in
Fig.~\ref{a&b}. If $Y$ is of type A we have $e(D)=4$, otherwise we have
$e(D)=6$, and the two discs of $F$ are completely symmetric.
\begin{figure}
\begin{center}
\includegraphics{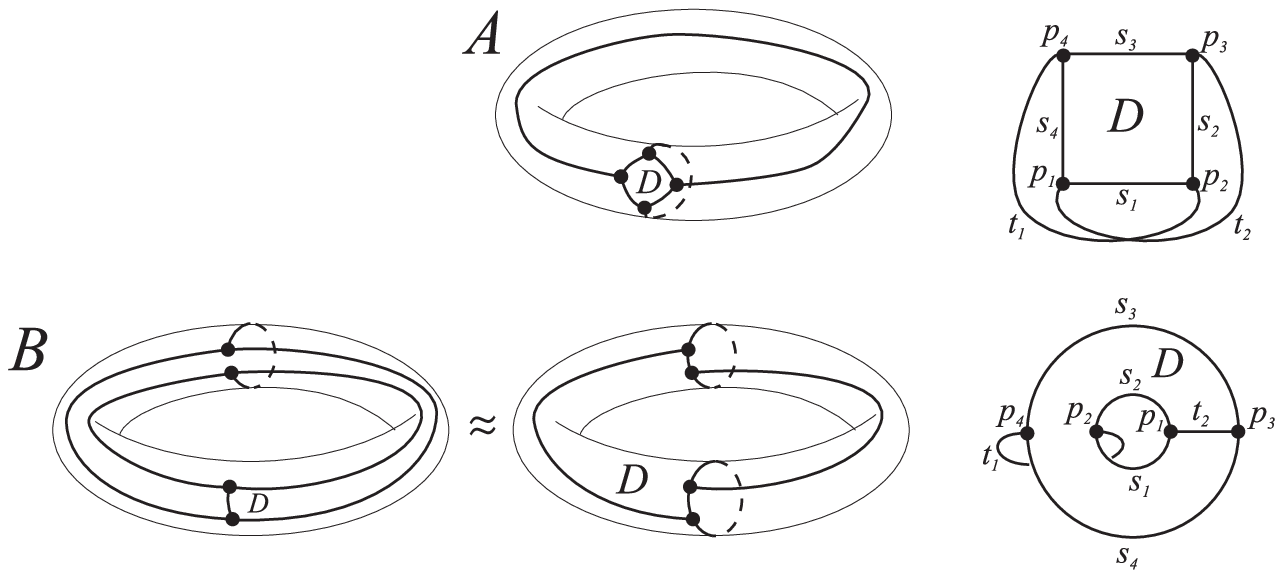}
\nota{Embeddings of type $A$ and $B$.} \label{a&b}
\end{center} \end{figure}
Figure~\ref{a&b} also contains notation used throughout the proof (note
that $s_1, \ldots, s_4$ are the edges in $\partial_1 D=\partial\Sigma_D$ both in case $A$
and in case $B$). Let $f_i$ be the face of $\Sigma_D \setminus S(P)$
incident to $s_i$. Moreover, let $g_j$ be the face of $P$ incident to
$t_j$. Since $D$ is good, we have $g_1, g_2 \not\subset\Sigma_{D}$.
Finally, let $e_i$ be the edge of $P$ which contains $p_i$.

\smallskip\noindent{\bf Fact 2.} \emph{Either the faces $f_1, f_2, f_3,
f_4, g_1, g_2$ are all distinct or $Y$ bounds a polyhedron of type
\ref{loop:item} or \ref{chain:item}.1.} \smallskip

Assuming that $Y$ does not bound a polyhedron of type \ref{loop:item} or
\ref{chain:item}.1, it follows that the segments $e_i\cap\Sigma_{D}$ for
$i=1, \ldots, 4$ are distinct. Then let $v_i \in V(P)$ be the endpoint of
$e_i\cap\Sigma_{D}$ not lying on $D$.

\smallskip\noindent{\bf Fact 3.} \emph{Up to symmetry we have $v_1=v_2$ in
case $A$ and either $v_1=v_2$ or $v_1=v_3$ in case $B$.} \smallskip

Let us now set $u=s_1$ in case $A$, and either $u=s_1$ or $u=t_1$ in case
$B$, depending on whether $v_1=v_2$ or $v_1=v_3$, so there are two edges
of $P\cup D$ which start at the endpoints of $u$ and both end at $v_1$.
These edges are $e'_1=e_1\cap\Sigma_D$ and $e'_m=e_m\cap\Sigma_D$, with
$m\in\{2,3\}$ depending on the case. Recall now that if two edges end at
the same vertex then one face incident to the first edge is also incident
to the second one. Since we are assuming that the $f_i$'s and $g_j$'s are
distinct, we deduce that $u\cup e'_1\cup e'_m$ bounds a disc of $P\cup D$,
which is a triangle, \emph{i.e.} $u\subset{\rm lk}(v_1)$. Following
Remark~\ref{J1:along} we can then perform a $J_1$-move to which
Lemma~\ref{J-lemma:2} and Lemma~\ref{J-lemma:3} apply. Denoting by $D'$
the disc corresponding to $D$ after the move, we have $e(D')\le e(D)$, and
equality can hold only if $Y$ is of type $B$.

\smallskip\noindent{\bf Fact 4.} \emph{If $e(D') < e(D)$ then $Y$ bounds
a polyhedron of type \ref{standard:item}.4 or \ref{standard:item}.5.}

\smallskip\noindent{\bf Fact 5.} \emph{If $e(D') = e(D)$ then $Y$ bounds
a polyhedron of type \ref{chain:item}.2.} \smallskip

This establishes the theorem. We now prove our assertions.

\smallskip\noindent{\bf Proof of fact 1.} By Theorem~\ref{shortloops} the
loop $\partial D$ is fake, and we can perform a move $J_{e(D)}$ as
explained in Subsection~\ref{J}. The result is a trace $Y'$ with 2 vertices
of a surface $F'$ isotopic to $F$. By Proposition~\ref{traces:2:easy} either
$F'$ is a non-separating torus, or $Y'$ is boundary-parallel, or we have
$Y'=\partial \calR(p)$ for some $p\in V(P)\setminus S(P)$. In the first
case, up to isotoping $F'$ back to $F$, getting an isotopic copy $P'$ of
$P$, we get a contradiction to our initial assumption. In the other cases
we have to see which polyhedra can result from an inverse $J_{e(D)}$ move
applied to $\theta_i\times[0,1]$ or to $\calR(p)$. It is now rather easy
to examine all possibilities and check the assertion.

\smallskip\noindent{\bf Proof of fact 2.} Of course no $f_i$ can be equal
to a $g_j$, because $f_i\subset\Sigma_D$ and $g_j\cap\Sigma_D=\emptyset$.
Let us first show that if two $f_i$'s coincide then $Y$ bounds a
polyhedron of type \ref{loop:item} or \ref{chain:item}.1. We refer to
Fig.~\ref{a&b} for the notation.

Two adjacent $f_i$'s cannot coincide because of
Corollary~\ref{no:embedded:faces}. Up to symmetry, the only cases we are
left to deal with are $A$-$(f_1=f_3)$, $B$-$(f_1=f_3)$, and
$B$-$(f_2=f_3)$. In all cases we will show that $Y$ bounds a polyhedron of
type 2 or 4.1. The key point will be to exhibit two loops that must be
fake because of Theorem~\ref{shortloops}.

Case $A$-$(f_1=f_3)$ is examined in Fig.~\ref{fact2}-left: since $\alpha'$
and $\alpha''$ are fake, one sees quite easily that
$Y=\partial\calR(\gamma)$, where $\calR(\gamma)$ is a M\"obius strip with
two tongues (Fig.~\ref{fact2}-right).
\begin{figure} \begin{center}
\includegraphics{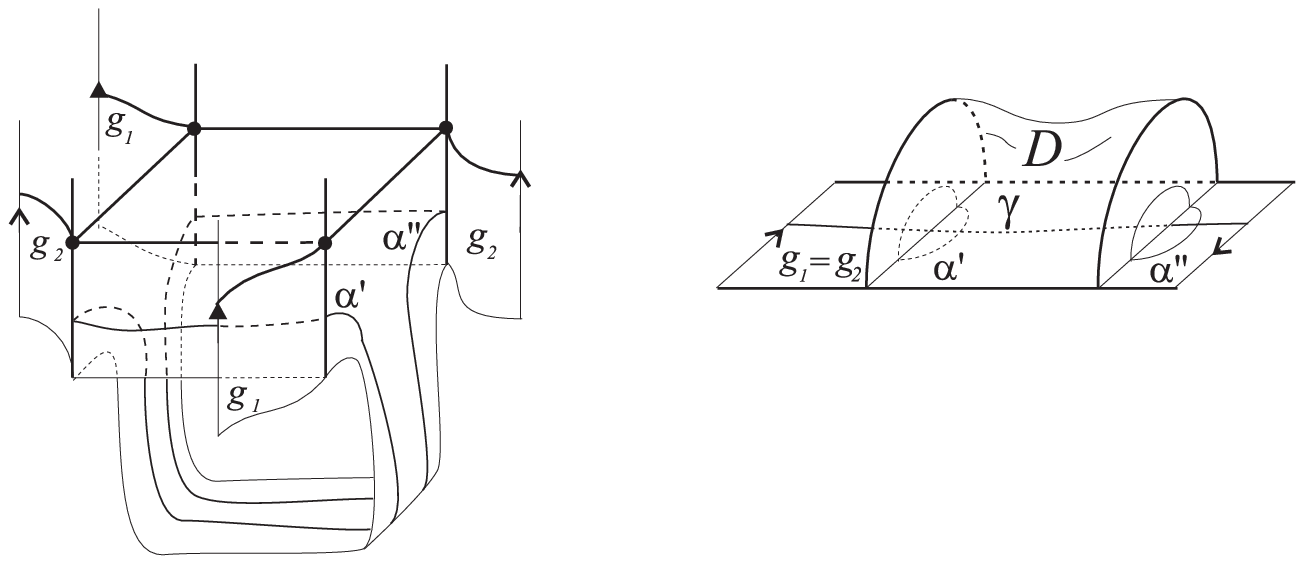}
\nota{Proof of fact 2: first case.} \label{fact2}
\end{center}
\end{figure}
Case $B$-$(f_1=f_3)$ is similar (Fig.~\ref{fact2b}-left); we have
$Y=\partial\calR(\gamma)$, where $\calR(\gamma)$ is an annulus with two
tongues on opposite sides (Fig.~\ref{fact2b}-right).
\begin{figure} \begin{center}
\includegraphics{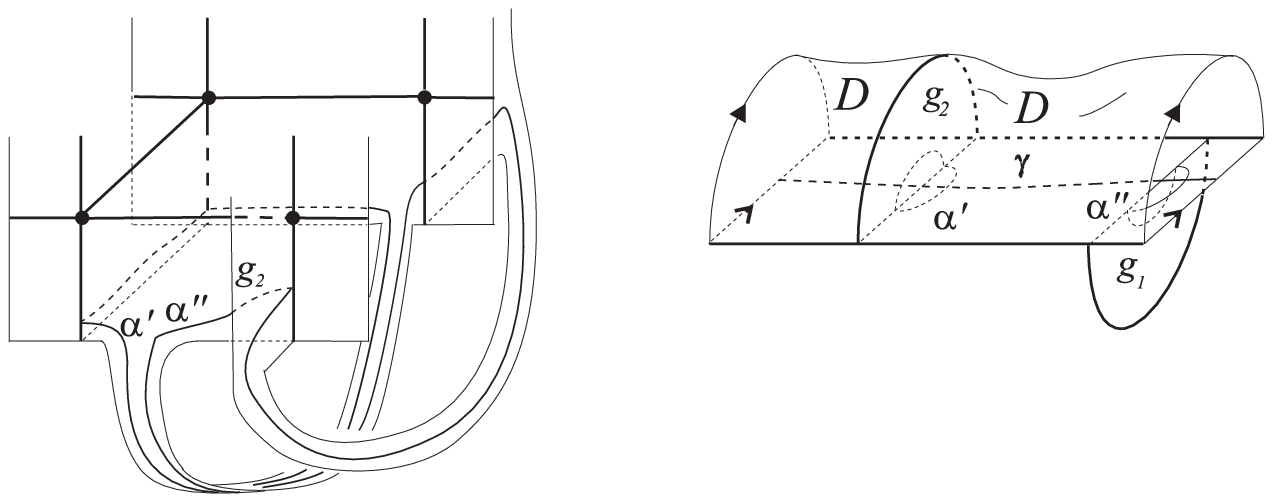}
\nota{Proof of fact 2: second case.} \label{fact2b}
\end{center}
\end{figure}
In case $B$-$(f_2=f_3)$ we consider the loops of Fig.~\ref{fact2c}-left.
Since $\alpha'$ and $\alpha''$ are fake we deduce that all the edges
$e_i\cap\Sigma_D$ end at the same vertex $v$, such that $s_2, s_3 \subset
{\rm lk}(v)$. We can then apply a move $J_1$ whose effect on $Y$ is shown in
Fig.~\ref{fact2c}-right.
\begin{figure} \begin{center}
\includegraphics{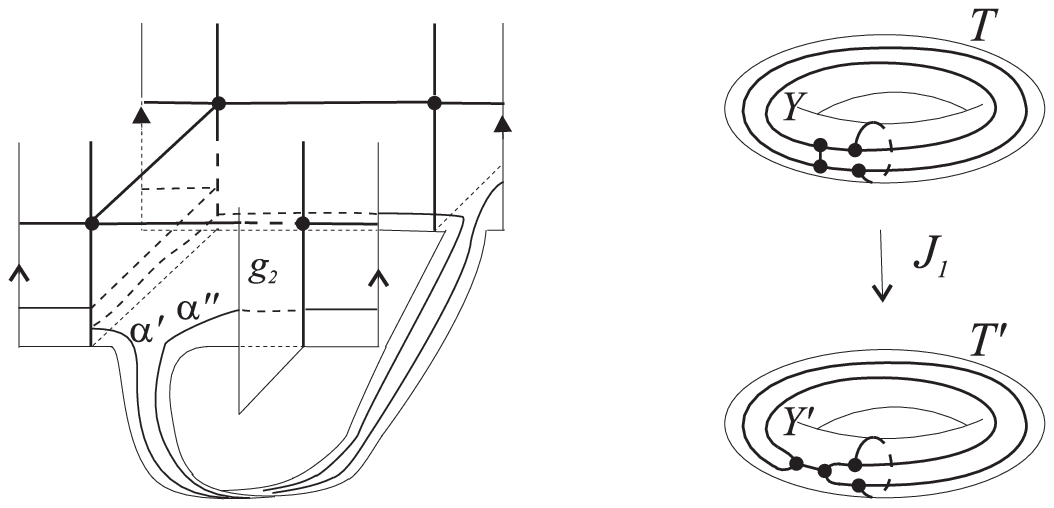}
\nota{Proof of fact 2: third case.} \label{fact2c}
\end{center} \end{figure}
The result is a trace $Y'$ which falls into case $A$-$(f_1=f_3)$. So
$Y'=\partial\calR(\gamma)$ with $\calR(\gamma)$ a M\"obius strip with two
tongues. Recalling that the inverse of a $J_1$-move is again a $J_1$-move,
we only need to consider which such moves can be applied to $\calR(\gamma)$. The move
is determined by the edge of $\partial\calR(\gamma)$ which disappears during the
move: of the 6 edges of $\partial\calR(\gamma)$, 4 lead to a situation in which
$e(D)=3$, so we exclude them. The other 2 edges are actually symmetric,
and the result is of type 4.1.

To conclude the proof of Fact 2 we must show that if the $f_i$'s are
distinct then $g_1\neq g_2$. If $Y$ is of type $B$ then $g_1$ has a
certain component of $M\setminus(P\cup D)$ on both sides, and $g_2$ has
the other one, so $g_1\neq g_2$. Assume in case $A$ that $g_1=g_2$.
Referring to Fig.~\ref{a&b} let $q_j$ be the midpoint of $t_j$, and join
$q_1$ to $q_2$ by an arc $\lambda$ in $g_1=g_2$. There are 4 distinct arcs
$\lambda_1, \ldots, \lambda_4 \subset Y$ having endpoints $q_1$ and $q_2$
and intersecting $S(P)$ twice. For two of them the polyhedron
$\calR(\lambda \cup \lambda_i)$ is an annulus with $2$ tongues on the same
side. Then some $\lambda\cup\lambda_i$ is fake, which is in contrast with
the fact that the $f_i$'s are distinct.

\smallskip\noindent{\bf Proof of fact 3.} We start with case $A$. Assume
that $v_1\neq v_2$, and put $P'=(P\cup D)\setminus f_1$. If $f_1$ is
incident to $x$ different vertices of $P$ then $\#V(P') = \#V(P)+4-2-x$.
Since $P$ is minimal we have $x\le 2$. On the other hand ${f}_1$ is
incident to $v_1$ and $v_2$, so $x=2$. Now Fig.~\ref{fact3d}
\begin{figure} \begin{center}
\includegraphics{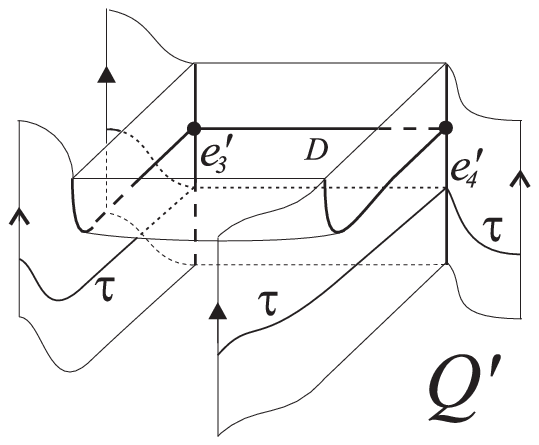}
\nota{A triod $\tau \subset P'$.}
\label{fact3d}
\end{center} \end{figure}
shows a triod $\tau$ in $P'$, trace of a torus parallel to $F$. By
Proposition~\ref{traces:2:easy}, either $F$ is non-separating or $\tau$ is
boundary-parallel. In the first case we get a contradiction to the initial
assumption. In the second case we deduce that $f_1$ is incident to $v_3$
and $v_4$, so $\{v_3,v_4\}\subset\{v_1,v_2\}$. So either $v_3=v_4$, or
$v_4=v_1$, or $v_3=v_1\neq v_4=v_2$. In all cases but the last one the
conclusion is the desired one up to symmetry. Concentrating on the last
case, we note that $f_1\cup\ldots\cup f_4$ is a surface near $v_1$ and
$v_2$, and that the $f_i$'s and $g_j$'s are all distinct. From these facts
it is not hard to deduce that $v_1$ and $v_2$ appear as in
Fig.~\ref{fact3}.
\begin{figure} \begin{center}
\includegraphics{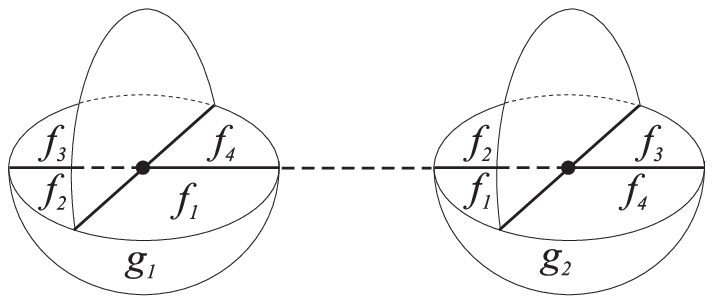}
\nota{The vertices $v_1$ and $v_2$.}
\label{fact3}
\end{center} \end{figure}
The figure readily implies that $f_2=f_4$: a contradiction.

The proof in case $B$ is similar, except that $D$ cannot be used directly:
a perturbed version $D'$ as in Fig.~\ref{fact3b}-left
\begin{figure} \begin{center}
\includegraphics{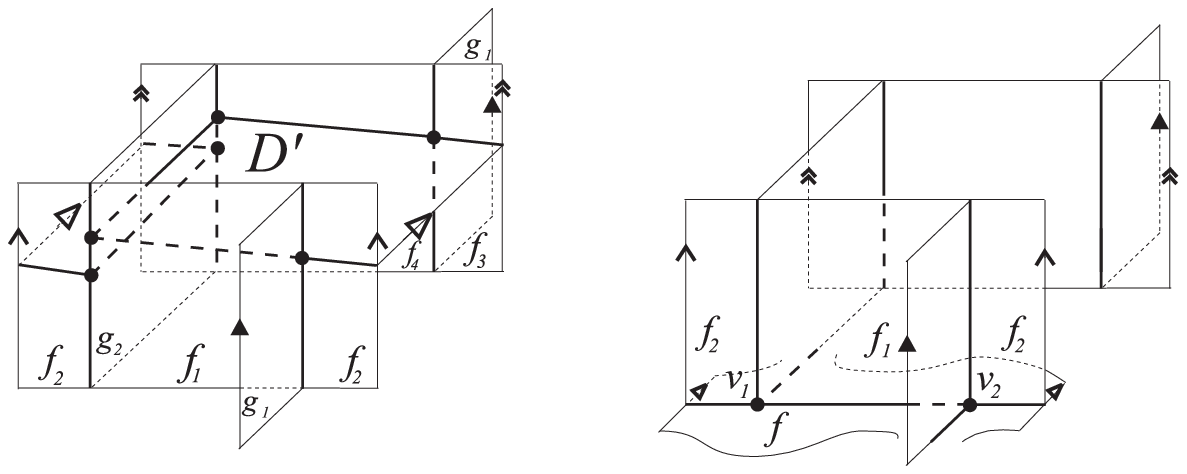}
\nota{The perturbed disc $D'$ (left), and a length-1 loop bounding an
external disc (right).} \label{fact3b}
\end{center} \end{figure}
must be employed. We are again supposing here that $v_1\neq v_2$, so
${f}_1$ is incident to $x\geq 2$ vertices of $P$, but now $\#V(P') =
\#V(P)+6-3-x$. Since $P$ is minimal we have $x\geq 3$, so $x\in\{2,3\}$.
We first claim that we can suppose $x=3$ up to symmetry. By contradiction,
assume that both ${f}_1$ and ${f}_2$ are incident to exactly 2 vertices.
We deduce that the situation is as in Fig.~\ref{fact3b}-right, where we
also show a face $f$ incident twice to an edge, which is absurd by
Corollary~\ref{no:embedded:faces}.

Our claim that $x=3$ up to symmetry is proved, so $\#V(P')=\#V(P)$ and
$P'$ is minimal too. A figure very similar to Fig.~\ref{fact3d} shows that
a triod must exist in $P'$, and allows to conclude as above that either
$F$ is separating or $f_1$ is incident to $v_3$ and $v_4$. So either
$v_3=v_4$, which gives the desired conclusion up to symmetry, or
$\{v_3,v_4\}\cap\{v_1,v_2\}\neq\emptyset$ (recall that ${f}_1$ is incident
to exactly 3 vertices). If $v_3=v_1$ or $v_4=v_2$ we get the desired
conclusion. Otherwise we can assume up to symmetry that $v_1=v_4$. So
$e_1$ and $e_4$ have a common vertex in $P$, which implies that there is a
face incident to both. But $e_1$ is adjacent to $f_1,f_2,g_1$ and $e_4$ is
adjacent to $f_3,f_4,g_2$, and the $f_i$'s and $g_j$'s are distinct, so we
get a contradiction.

\smallskip\noindent{\bf Proof of fact 4.} If $Y$ is of type $A$, then
$e(D')=3$, so by Fact 1 (and its proof) $Y'$ bounds a polyhedron $Q$ of
type \ref{standard:item}.3 or of type \ref{vert:nhb:item}, but the latter
is impossible because $Y'$ is the trace of a torus. We only need to
consider which $J_1$-moves can be applied to a $Q$ of type
\ref{standard:item}.3. By Lemma~\ref{J-lemma:3} the move actually takes
place towards the exterior of $Q$ (\emph{i.e.}~its result contains 2
vertices of $P$). The move is determined by the edge of $\partial Q$
which disappears during the move: of the 6 edges in $\partial Q$, 3 lead
to a situation in which $e(D)=2$, so we exclude them. The other 3 edges
are actually symmetric, and the result is one of the polyhedra of type
\ref{standard:item}.4.

If $Y$ is of type $B$, then $u$ must be an edge in $\partial_2D$
(otherwise $e(D')=e(D)$), so $Y'$ is of type $A$. Moreover $Y'$ is the
trace of a torus. Combining Fact 2 and the part of Fact 4 already
established we see that $Y'=\partial Q$ with $Q$ either of type
\ref{standard:item}.4 or a M\"obius strip with two tongues (type 2).
However, if we denote by $f_i'$ the faces of $\Sigma_{D'}$ incident to
$D'$, by Lemma~\ref{J-lemma:3} we have $f'_i\subset f_i$ up to
permutation, so the $f_i'$'s are distinct. This shows that type 2 is
impossible, and again we are left to analyze what can we get from a $Q$
of type \ref{standard:item}.4 by a move $J_1$ which takes place towards
the exterior. Of the 6 edges of $\partial Q$, 4 lead to a situation in
which $e(D)=3$, so we exclude them. The other 2 edges are actually
symmetric, and the result is type \ref{standard:item}.5.

\smallskip\noindent{\bf Proof of fact 5.} The first step of our proof is
the extension of the move $Y\to Y'$ to a flow $Y\to Y'\to Y''\to\cdots\to
Y^{(k)}$ of $J_1$-moves. As mentioned in the proof of Fact 4 we must have
$u\subset\partial_1 D$ in this case, so we assume up to symmetry that
$u=s_1\subset\partial f_1$, and we note that
Remarks~\ref{J-lemma:1}-\ref{ultimo:ritocco?} and Lemma~\ref{J-lemma:3}
apply. The situation is described in Fig.~\ref{fact5b}.
\begin{figure} \begin{center}
\includegraphics{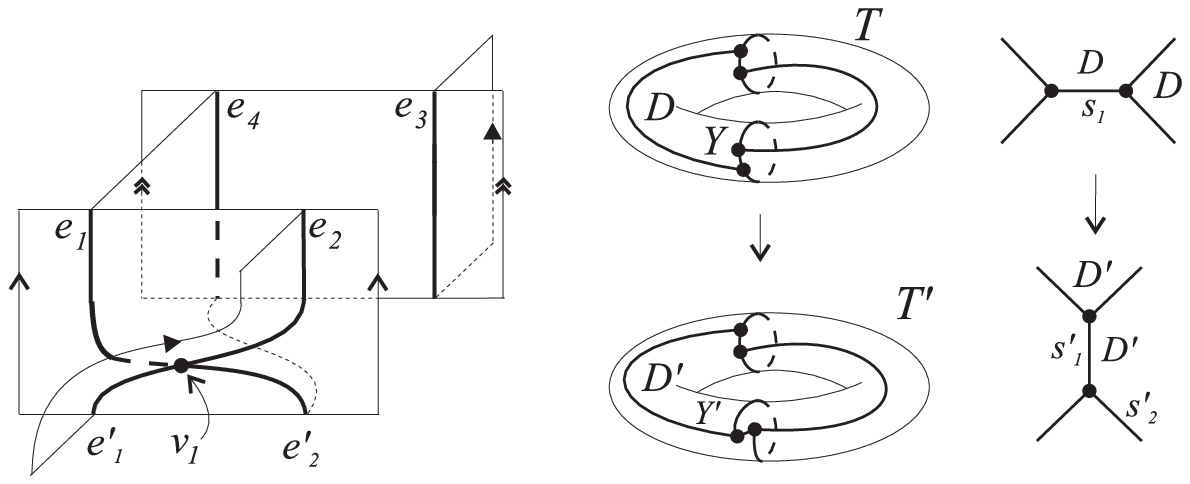}
\nota{Proof of Fact 5.} \label{fact5b}
\end{center}
\end{figure}
One easily sees that the faces of $\Sigma_{D'}$ incident to $\partial_1
D'$ are $f_3,f_4$ and two new ones (one of which is contained in $f_2$),
which we denote by $f_1',f_2'$. If $\{f_1',f_2',f_3,f_4\}$ are not
distinct, the flow is reduced to $Y\to Y'$, and we move to the next step.
Otherwise let $v_1',v_2'$ be the ends of $e_1',e_2'$ (see
Fig.~\ref{fact5b}-left). If $v_1'\neq v_2'$ then again the flow is reduced
to $Y\to Y'$. Assume on the contrary that $v_1'=v_2'$, and consider
Fig.~\ref{fact5b}-right. Then either $s_1'$ or $s_2'$ is contained in
${\rm lk}(v')$, but certainly $s_1'$ is not, for otherwise $P$ would
contain an embedded face with two vertices, which is absurd by
Corollary~\ref{no:embedded:faces}. Setting $u'=s_2'$, we are now in a
position to apply a move $J_1$ along the triangle determined by $v_1'$ and
$u'$, getting from $Y'$ to $Y''$. We proceed in a similar way and note
that the process must come to an end because $\Sigma_{D^{(i)}}$ contains
one vertex less than $\Sigma_{D^{(i-1)}}$ by Lemma~\ref{J-lemma:3}.

Our second step is to understand the final stage $Y^{(k)}$ of our flow. By
construction either $\{f_1^{(k)},f_2^{(k)},f_3,f_4\}$ are not distinct or
$v_1^{(k)}\neq v_2^{(k)}$. In the first case, since at each step only 1
face not contained in the previous one is inserted (and 1 is deleted),
precisely 3 of $\{f_1^{(k)},f_2^{(k)},f_3,f_4\}$ are distinct. We know by
Fact 2 (and its proof) that $Y^{(k)}$ (which is of type $B$) bounds a
polyhedron $Q$ which is either an annulus with 2 tongues on opposite
sides, or of type 4.1. The first case is excluded by what just said about
the $f_i$'s. By Lemma~\ref{J-lemma:3}, $Y$ bounds
$Q\cup_{Y^{(k)}}[Y,Y^{(k)}]$. Since at each step of the construction of
our flow the choice of move $J_1$ was forced, the polyhedron $[Y,Y^{(k)}]$
is defined unambiguously (it depends on $k$ only). We only need to explain
which edge of $\partial Q$ determines the $J_1$-move which glues $Q$ to
$[Y,Y^{(k)}]$. Of the 6 edges, 2 lead to a trace of type $A$, 2 give rise
to an embedded face with 2 edges (excluded by
Corollary~\ref{no:embedded:faces}) and the other 2 are symmetric, so
$Q\cup_{Y^{(k)}}[Y,Y^{(k)}]$ also depends on $k$ only. It is now a routine
matter to check that indeed $Q\cup_{Y^{(k)}}[Y,Y^{(k)}]$ is the polyhedron
of type 4.2 with $k$ vertices.

Having understood the case where $\{f_1^{(k)},f_2^{(k)},f_3,f_4\}$ are not
distinct, we assume that they are. The rest of the proof is devoted to
showing that it is actually impossible that $v_1^{(k)}\neq v_2^{(k)}$. Let
us first assume that $v_3\neq v_4$. By Fact 3 we then have $v_1=v_3$ up to
symmetry, and we can apply a move $J_1$ which reduces $e(D)$. Fact 4 shows
that $Y^{(k)}$ bounds a polyhedron $Q$ of type \ref{standard:item}.4 or
\ref{standard:item}.5, but $\partial Q$ is of type $B$, so it must be of
type \ref{standard:item}.5. Once again we must analyze the possible
results of a move $J_1$, towards the exterior of a $Q$ of type
\ref{standard:item}.5. Of the 6 edges of $\partial Q$, 2 lead to a trace
of type $A$, and therefore are excluded. The 4 other edges come in 2
symmetric pairs. For one type, the result of the move $J_1$ contains an
embedded face with 3 vertices, which is absurd by Corollary~\ref{no:embedded:faces}.
\begin{figure} \begin{center}
\includegraphics{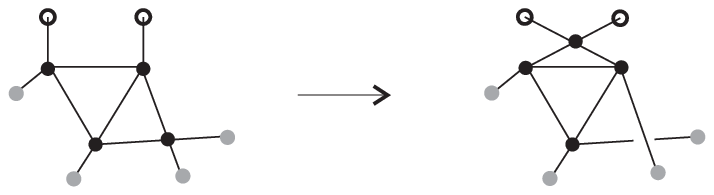}
\nota{Proof of Fact 5 continued.} \label{fact5} \end{center}
\end{figure}
For the other type, the result contains an embedded face with 4 vertices.
We can then apply a disc-replacement move as in Fig.~\ref{badsquare},
getting a new minimal skeleton $P'$ of $(M,X)$.
The evolution of the singular set is shown in Fig.~\ref{fact5},
where the two white dots lie on some $\theta_i$, the black dots are
vertices, and the gray dots lie on $Y$. Since the edges leaving
$\theta_i$ end at the same vertex, a $J_1$-move transforms
$\theta_i$ into a triod which is not boundary-parallel.
This contradicts Proposition~\ref{traces:2:easy}.

We are left to deal with the case where $\{f_1^{(k)},f_2^{(k)},f_3,f_4\}$
are distinct, $v_1^{(k)}\neq v_2^{(k)}$, and $v_3= v_4$. In this case we
can perform a $J_1$-move along either $s_3$ or $s_4$, and we can proceed
just as above, constructing a flow $Y^{(k)}\to Y^{(k+1)}\to \dots\to
Y^{(k+h)}$. During this process the faces $f_1^{(k)},f_2^{(k)}$, and the
vertices $v_1^{(k)},v_2^{(k)}$ remain unaffected, while $f_3,f_4,v_3=v_4$
get transformed into $f_3^{(h)},f_4^{(h)},v_3^{(h)},v_4^{(h)}$. As above,
we have at the end of the sequence either that
$\{f_1^{(k)},f_2^{(k)},f_3^{(h)},f_4^{(h)}\}$ are not distinct or that
$v_3^{(h)}\neq v_4^{(h)}$. In the first case, Fact 2 implies that
$Y^{(k+h)}$ bounds a polyhedron of type 2 or 4.1. Such a polyhedron has at
most 1 vertex, but $\Sigma_{D^{(k+h)}}$ contains at least $v_1^{(k)}\neq
v_2^{(k)}$, and we get a contradiction. In the second case we are
precisely in the situation $v_3\neq v_4$ previously considered, and again
we get a contradiction. \end{proof}

\subsection{Conclusion of proofs}\label{real:properties}

If $Y$ is a trace in $P$, we denote by $P_Y$ the polyhedron $P\setminus
\calR(Y)$.

\medskip\noindent \emph{Proof of Theorem~\ref{minimal:are:standard}.} Let
$P$ be a minimal skeleton of $(M,X)$.
By Corollary~\ref{zero:bricks} we have $c(M,X)>0$, so $P$ is standard.
Suppose a face $f$ of $P$ is
incident to $\partial P$ in at least two distinct edges $e \subset
\theta_i$ and $e' \subset \theta_{i'}$. We note that $i\neq i'$ by
Lemma~\ref{faces:distinct}, and choose an arc $\alpha$ in $f$ having one
end on $e$ and one on $e'$. Then $Y=\partial\calR(\theta_i \cup
\theta_{i'} \cup \alpha)$ is a trace with 4 vertices of a surface $F$.
Moreover $P_Y = P_1 \sqcup P_2$ is disconnected, so $F$ separates $M$ and
hence it is orientable. Let $P_2$ be the component containing $\alpha$.

The graph $Y$ is of type B (see Fig.~\ref{abgraphs}) and $P_2$ has 3
boundary components (namely, $\theta_i$, $\theta_{i'}$, and $Y$). Now
either $P_1$ or $P_2$ is of one of the types listed by
Theorem~\ref{traces:4}, but no such type has 3 boundary components, so
$P_1$ must be of one such type. The only polyhedra among those listed in
Theorem~\ref{traces:4} having at least one vertex and boundary of type B
are those of type 3.5 (Fig.~\ref{type3}) and 4 (Fig.~\ref{type4}).
If $P_1$ is of type 3.5 then
$P$ is the skeleton of $B_4$, and if $P_1$ is of type 4 with 1 vertex then
it is the skeleton of $B_3$. Otherwise $P_1$ is of type 4 with $k\geq 2$
vertices, and the two edges of $S(P)$ adjacent to $\theta_i$ have a
common endpoint. It easily follows that via a $J_1$-move we can transform
$\theta_i$ into a triod which is not boundary-parallel and is the trace of
a separating torus. This contradicts Proposition~\ref{traces:2:easy}.


\vspace{-.85cm}
\begin{flushright} $\square$ \end{flushright}

\medskip\noindent \emph{Proof of Theorem~\ref{disjoint:length:1:loops}.}
Set $\calL=\{\gamma_1,\ldots,\gamma_n\}$, where $\gamma_i$ is the core of
the M\"obius strip with one tongue attached to $\theta_i\subset\partial
P$. By Theorem~\ref{minimal:are:standard}, even if we modify each $\gamma_i$
within its isotopy class, the $\gamma_i$'s stay disjoint. Moreover, each
$\calR(\gamma_i)$ is a M\"obius strip with one tongue. Therefore it is enough
to show that $\calL$ is a set of representatives of length-1 loops in $Q$.
If not, there is a length-1 loop $\gamma$ not isotopic to any $\gamma_i$.

If $\gamma$ is disjoint from all $\gamma_i$'s, then $\gamma\subset
P$, so a face of $P$ is doubly incident to some edge, and we get a
contradiction to Corollary~\ref{no:embedded:faces}. If $\gamma$ meets
some $\gamma_i$ then, by Theorem~\ref{minimal:are:standard}, it meets
only one, and we can assume that $\gamma\cap\gamma_i$ is one point
away from $S(Q)$. Set $R=\calR_Q(\gamma\cup\gamma_i)$. We need now to
distinguish two cases, depending on whether $\calR_Q(\gamma)$ is a
M\"obius strip or an annulus with one tongue. In the first case there
exists a curve $\alpha$ contained in $\partial R$, and therefore in
$P$, such that $l(\alpha)=2$ and $\alpha$ bounds an external disc
($\alpha$ is homologous to $\gamma+\gamma_i$ in $R$). By
Theorem~\ref{shortloops} $\alpha$ is fake, and it easily follows that
$\gamma$ is isotopic to $\gamma_i$.

Assume now that $\calR_Q(\gamma)$ is an annulus with one tongue. Note
that $\partial R\subset P$ is a trace with 4 vertices of a
separating, and hence orientable, surface $F$. Moreover $\partial R$
is of type $A$, so, by Theorem~\ref{traces:4}, $\partial R$ bounds in
$P$ a polyhedron $S$ of type 1.1, 3.3, 3.4, or 2 based on a M\"obius
strip. But $R\cap P$ is not of such a type, so the rest of $P$ is,
hence $\#V(P)\le 1$. But $\calB^1_{\le 1}=\{B_0,\ldots,B_3\}$, and we
are done.


\vspace{-.85cm}
\begin{flushright} $\square$ \end{flushright}

\medskip\noindent
Before proving Theorem~\ref{restrictions:on:graph} we establish a general fact.

\begin{lemma} \label{orthogonal} Let $Q$ be a filling of a minimal
skeleton $P$ of a brick. Let $\{e_1, \ldots, e_{2m}\}$ be a set of edges
which disconnects $S(Q)$ in two components. Then there is a trace $Y$
contained in $P$ which has $2m$ vertices $p_i\in e_i$ for $i=1,\ldots,2m$,
and $Y$ is the trace of an orientable separating surface. \end{lemma}

\begin{proof} Take points $p_i\in e_i$; we have $S(Q)\setminus\{p_i\} = K_1 \sqcup
K_2$. Let $f$ be a face of $Q$ incident to some $e_i$. The gluing path of
$\partial f$ to $S(Q)$ can be split into arcs $s_1, \ldots s_{2\nu}$,
meeting at points $q_1,\ldots,q_{2\nu}$, where $s_{2j+1} \subset K_1$ and
$s_{2j} \subset K_2$ for all $j$, and each $q_k$ is glued to one
$p_{\beta(k)}$. The map $\beta$ is not necessarily injective, since $f$
can be multiply incident to an edge $e_i$. We can give the points $q_k$
alternating (red and black) colors.

Since $P = Q\setminus\calR(\calL(Q))$ is super-standard, $f$ can intersect
at most one loop $\gamma$ among those in $\calL(Q)$. Now take $\nu$
pairwise disjoint segments $\lambda_1, \ldots, \lambda_\nu$, properly
embedded in $f$, such that $\cup_{j=1}^\nu \partial \lambda_j =
\cup_{k=1}^{2\nu} q_k$. We can ask the $\lambda_j$'s to be disjoint from
$\gamma$, since the points on $\partial f$ are separated into two even
subsets by $\gamma$. It is easy to see that the two endpoints of each
$\lambda_j$ automatically have distinct colours. If we do this for each
face $f$ incident to some $e_i$, the union of all the chosen segments is a
trace $Y$ disjoint from $\calL(Q)$ and hence contained in $P$.

We claim that $Y$ has a product regular neighbourhood in $P$: take
for $i=1,\ldots, 2m$ a vector $v_i$ at $p_i$, tangent to $e_i$ and
directed towards $K_2$. Each segment of $Y$ is a $\lambda_j$,
properly embedded in a face $f$ such that $\partial \lambda_j$
consists of points with distinct colors. It follows that the vectors
at the ends of $\lambda_j$ extend along $\lambda_j$ to a
non-vanishing field tangent to $f$. The existence of such a field on
$Y$ easily implies that $F$ is orientable and that $F$ cuts $M$ into
two components.\end{proof}

\medskip\noindent \emph{Proof of Theorem~\ref{restrictions:on:graph}.}
Suppose $S(Q)$ contains a pair $\{e_0, e_1\}$ of separating edges. By
Lemma~\ref{orthogonal} there is a trace $Y$ of a separating (and hence
orientable) surface $F$ with two vertices, intersecting both $e_0$ and $e_1$.
Proposition~\ref{traces:2:easy} applies, and possibility (1) is ruled
out because $F$ separates. Both other possibilities imply that
the vertices of $Y$ lie on the same edge of $Q$, but $e_0\ne e_1$ by assumption.

Suppose $S(Q)$ contains a separating quadruple $\{e_0, e_1, e_2,
e_3\}$ of edges. By Lem\-ma~\ref{orthogonal} there is a trace $Y$ of
a separating (and hence orientable) surface with 4 vertices
intersecting them. If $Y$ is connected then Theorem~\ref{traces:4}
applies, and we are done because the singular sets of polyhedra of
types 1-4 indeed are as shown in Fig.~\ref{2bridge}. If $Y=Y_0\sqcup
Y_1$, then $Y_0$ is a trace with two vertices to which
Proposition~\ref{traces:2:easy} applies. Now possibility (1) is ruled
out either because every torus in $M$ is separating or by
Theorem~\ref{traces:2:difficult}, and as above the other two
possibilities lead to a contradiction.

\

\vspace{-.85cm}
\begin{flushright} $\square$ \end{flushright}

\section{Bricks and skeleta up to complexity 9}\label{concluding}

We provide in this section a complete description of the bricks in
$\calB_n$ for $n\le 9$ anticipated in Subsection~\ref{results}. Recall
that $\calB_n$ was split as $\calB^0_n\sqcup\calB^1_n$, where $\calB^0_n$
consists of the elements of $\calB_n$ without boundary. We describe now
$\calB^1_{\le 9}$, postponing $\calB^0_{\le 9}$ for a moment, because to
discuss it we will first need to introduce a new move on skeleta.

\begin{figure} \begin{center}
\includegraphics[width = 12 cm]{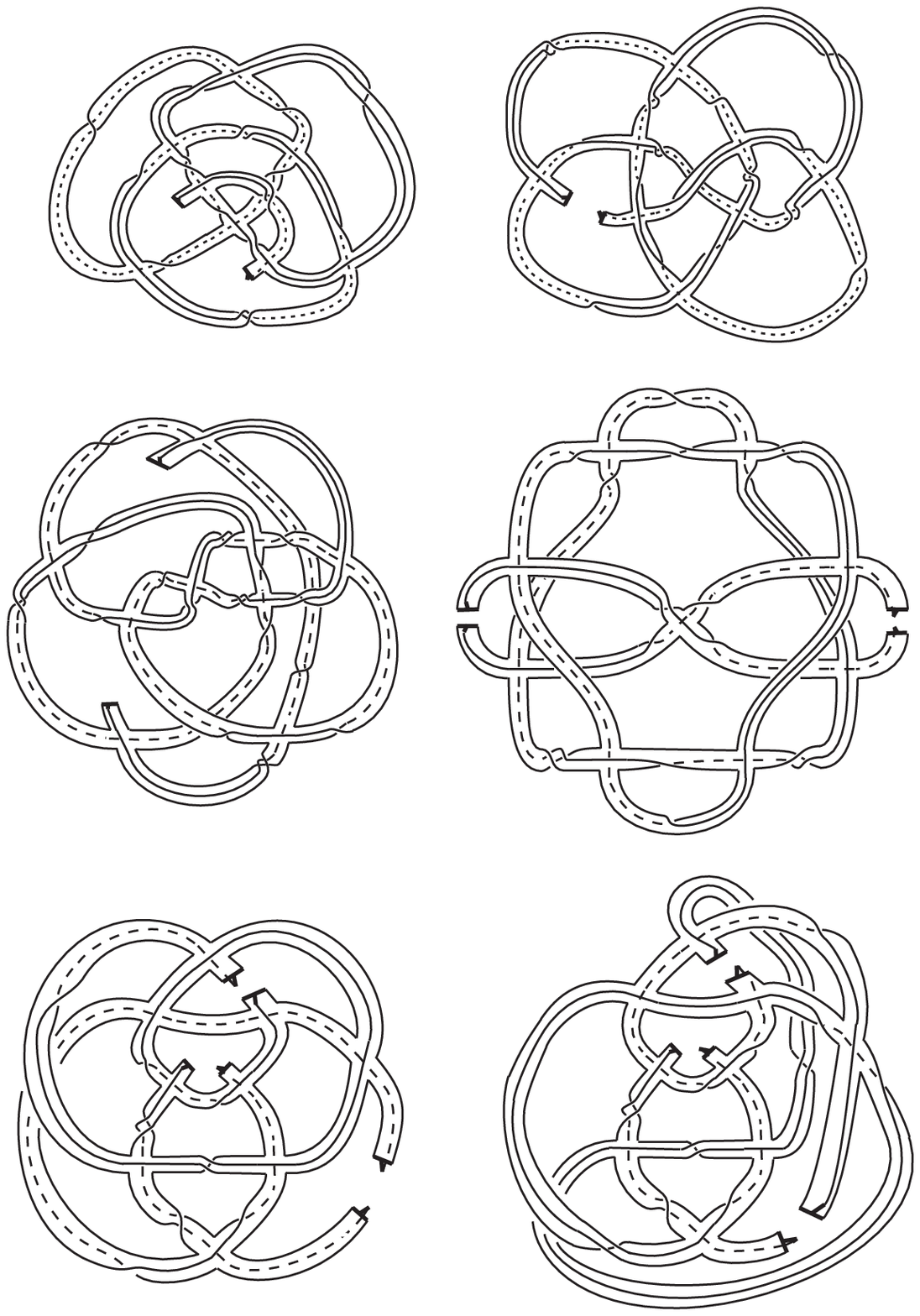}
\nota{Minimal skeleta for $B_5,\ldots B_{10}$.} \label{bigfig}
\end{center}
\end{figure}

Our computations show that the set $\calB_{\le 9}^1$ consists of $11$
bricks $B_0, \ldots, B_{10}$. Moreover, for $i\le 9$ there is a
unique minimal skeleton of $B_i$, while for $i=10$ there are two.
Minimal skeleta for
 $B_0,\ldots,B_4$ were shown in
Figg.~\ref{smallbr}
 and~\ref{smallbr2}, and for $B_5,\ldots,B_{10}$ they are
now shown in
 Fig.~\ref{bigfig}.
Using Remark~\ref{pictures}, in this figure we only draw
$\partial\calR_P(S(P))$, and we use a thicker line for the $Y$-shaped
portions of $\calR_P(S(P))$ lying on $\partial P$. Each component of
$\partial B_i$ contains two such $Y$'s (shown close to each other
when $\partial B_i$ has more than one component).

Having described $B_0,\ldots,B_{10}$, we can now
prove Proposition~\ref{not:sharp2}.

\medskip\noindent \emph{Proof of Proposition~\ref{not:sharp2}.}
Suppose $(M,\emptyset)$ is a sharp assembling of $B_i$ with $i\ge 6$
and some $B_2$'s and $B_3$'s. Since $c(M)\le 9$, one $B_3$ can occur
if $i=6$ only. Let $P_i$ be the minimal skeleton of $B_i$ shown in
Fig.~\ref{bigfig}. A minimal skeleton $P$ for $M$ is then a filling
of $P_i$, possibly after glueing one copy of the minimal skeleton of
$B_3$ if $i=6$. If we check all the polyhedra which can be built in
this way, we see that many of them contain embedded faces with no
more than 3 vertices, which contradicts Theorem~\ref{shortloops}.
Only 16 of them do not contain
 such a face.
Now 9 of these 16 are shown to be non-minimal by checking that small
faces appear after suitable
disc-replacement moves. The 7 polyhedra left out are skeleta
 of the 4
mentioned hyperbolic manifolds (there are some duplicates).


\vspace{-.85cm}
\begin{flushright} $\square$ \end{flushright}

\subsection{Twists} We introduce here a notion needed below to describe
$\calB^0_{\le 9}$. Let $P$ be a quasi-standard skeleton of a closed
manifold $(M,\emptyset)$, and let $\gamma$ be a length-$2$ loop in
$P$ such that $\calR(\gamma)$ is an annulus with $2$ tongues. For
$k\ge 1$ let $W_k$ be the polyhedron of type 4 with $k$ vertices
(Fig.~\ref{type4}). The boundaries $\partial \calR(\gamma)$ and
$\partial W_k$ are homeomorphic (of type $B$). We can then choose a
homeomorphism $\psi:\partial W_k \to \partial \calR(\gamma)$ and form
a polyhedron $P_k = {P\setminus \calR(\gamma)}\cup_{\psi} W_k$. Note
now that $W_k$ naturally sits in a solid torus $H$, with $\partial
W_k = W_k\cap \partial H$.

\begin{prop}\label{twist:prop} The homeomorphism $\psi:\partial W_k \to
\partial \calR(\gamma)$ can be chosen so that it extends to a
homeomorphism $\Psi:\partial H \to \partial \calR_M(\gamma)$. For
these choices $P_k$ is a skeleton of the Dehn surgered manifold $M_k =
{M\setminus\calR_M(\gamma)}\cup_{\Psi}H$. \end{prop}

\begin{proof} The first assertion is easy and taken for granted. By
construction $P_k$ sits in $M_k$ and it is simple, so we only need to
show that $M_k\setminus P_k$ is an open $3$-ball. To this end we note
that $M\setminus{P\cup\calR_M(\gamma)}$ is a ball $B$. Moreover
$\partial H\setminus \partial W$ consists of two discs $D'$ and $D''$, and
$H\setminus(\partial H \cup W)$ consists of two balls $B'$ and $B''$, with
$\partial{B'}\cap \partial H = {D'}$ and
$\partial{B''}\cap\partial H = {D''}$. So $M_k\setminus P_k = B
\cup_{\Psi|_{D'}}B' \cup_{\Psi|_{D''}}B''$ is a ball.
\end{proof}

We say that $P_k$ is obtained from $P$ by a $k$-\emph{twist} along
$\gamma$, and we adopt the convention that making a $0$-twist means
leaving $P$ unaffected.

\subsection{Closed bricks up to complexity 9}

\begin{figure}
\begin{center}
\includegraphics{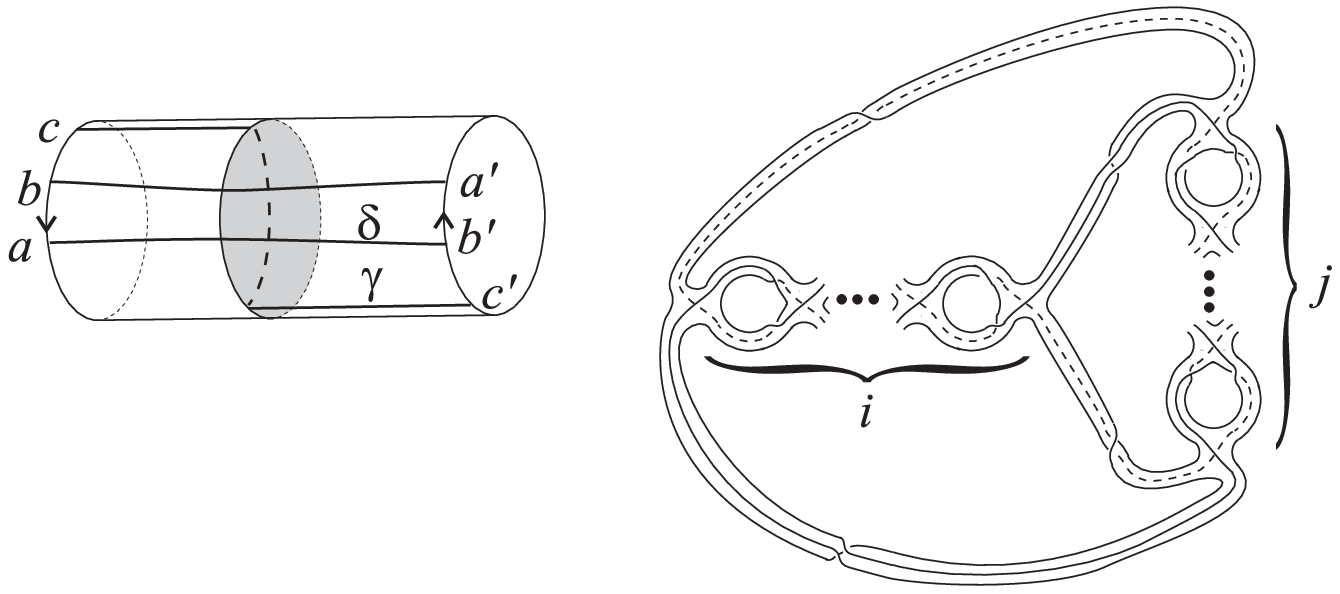}
\nota{The length-2 loops $\gamma$ and $\delta$ in $\tilde{C}_{0,0}$ and the
skeleton $\tilde{C}_{i,j}$.} \label{Cij}
\end{center} \end{figure}

Our computations show that the set $\calB_{\le 9}^0$ consists of $19$
bricks which belong to the union of two classes $\{C_{i,j}\}$ and
$\{E_k\}$. We describe here these manifolds and minimal skeleta
$\tilde{C}_{i,j}$ and $\tilde{E}_k$ of them. (As opposed to the case
of $\calB_{n\le 9}^1$, minimal skeleta are
often not unique in $\calB_{n\le 9}^0$.) The polyhedron $\tilde{C}_{0,0}$
of Fig.~\ref{Cij}-left is a
skeleton of $(S^2 \times S^1,\emptyset)$ and it contains $2$ length-2
loops $\gamma$ and $\delta$, shown in
Fig.~\ref{Cij}-left, such that ${S^2 \times S^1 \setminus
\calR_{S^2 \times S^1}(\gamma \cup \delta)} \cong (A, (2,1))$, where $A$
is the annulus. Both $\calR_{\tilde{C}_{0,0}}(\gamma)$ and
$\calR_{\tilde{C}_{0,0}}(\delta)$ are annuli with two tongues on different
sides. We can therefore perform an $i$-twist along $\gamma$ and a
$j$-twist along $\delta$. If we do this with appropriate gluing maps we
get the skeleton shown in Fig.~\ref{Cij}-right, which we denote by
$\tilde{C}_{i,j}$.
Using Proposition~\ref{twist:prop} it is not hard to check that
$\tilde{C}_{i,j}$ is a skeleton of the Seifert manifold $C_{i,j}=(S^2,
(2,1), (1+i,1), (1+j,1), (1,-1))$. We have $C_{i,j} = C_{j,i}$ for all
$i,j$.

\begin{figure}
\begin{center}
\includegraphics{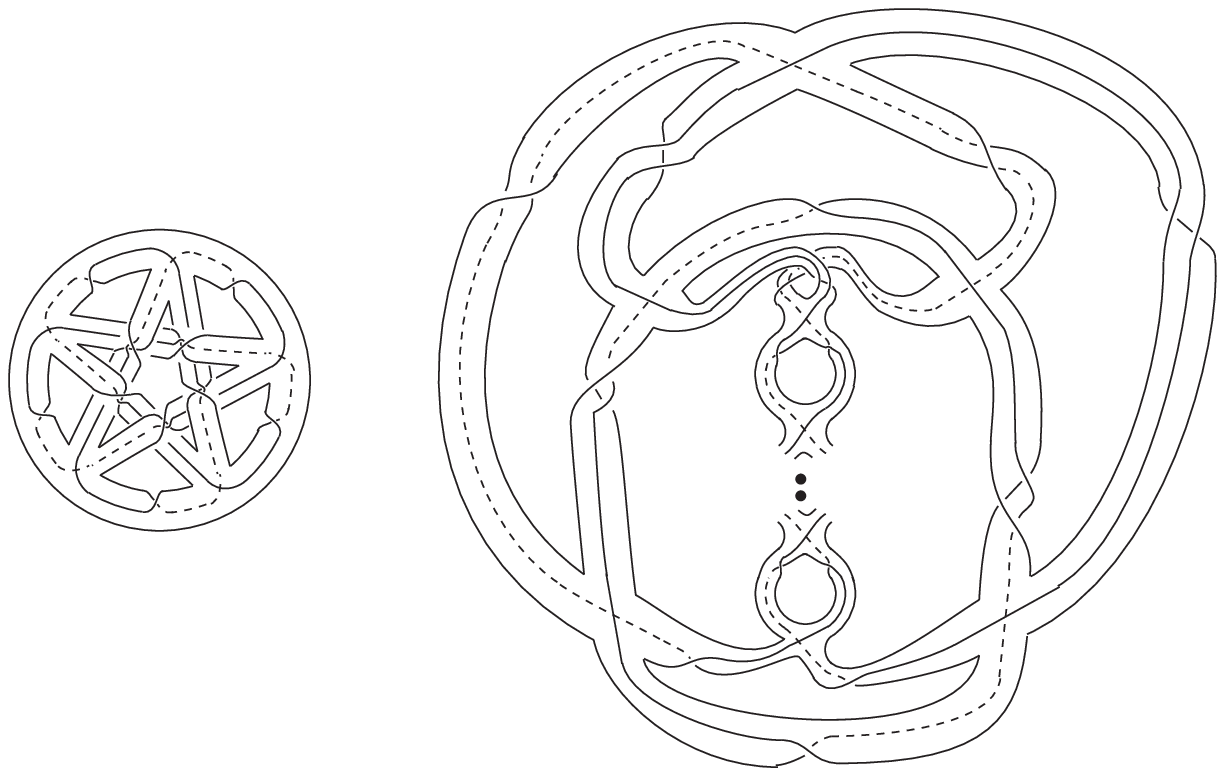}
\nota{The skeleta $\tilde{E}_0$ and $\tilde{E}_k$.} \label{Dk}
\end{center} \end{figure}

Poincar\'e's homology sphere $(S^2,(2,1), (3,1), (5,1), (1,-1))$ has a
unique minimal skeleton $\tilde{E}_0$ (Fig.~\ref{Dk}-left). For any pair
of non-adjacent edges of $S(\tilde{E}_0)$ there is a length-2 loop
$\gamma$ intersecting them, isotopic to the singular fiber $(5,1)$. Since
$\calR(\gamma)$ is an annulus with two tongues, we can perform a $k$-twist
along $\gamma$. If we do this with an appropriate gluing map we get the
skeleton shown in Fig.~\ref{Dk}-right, which we denote by $\tilde{E}_k$.
Each $\tilde{E}_k$ turns out to be a skeleton of the manifold
$E_k=(S^2,(2,1), (3,1), (5+k,1), (1,-1))$. It is worth mentioning here that
the
 minimal skeleton of the brick $B_5$ may be obtained from $\tilde{E}_0$ by
an
 operation similar to a $k$-twist along $\gamma$, except that the
polyhedron of type 3.5 (Fig.~\ref{type3}) is employed instead of $W_k$.

The set $\calB_{\le 9}^0$ consists of all manifolds $C_{i,j}$ and $E_k$
with $k \ge 0$ and $i\ge j \ge 1$ having at most $9$ vertices
(\emph{i.e.}~with $k\le 4$ and $i+j\le 9$), except the cases $k=1$ and $(i\ge
4,j=2)$. The skeleton $\tilde{E}_1$ is indeed minimal, but the associated
manifold is not a brick, since it lies in $\langle B_0\rangle_{\rm self}$.
This is coherent with the well-known fact that $(S^2,(2,1),(3,1),(6,1))$
fibers over $S^1$ with torus fiber. Each $\tilde{C}_{i,0}$ is minimal (for
$i\le 9$), but the corresponding manifold is contained in $\langle B_2,
B_3 \rangle_{\rm non-self}$. Each $\tilde{C}_{i,2}$ for $i\ge 4$ is not
minimal, since $\tilde{E}_{i-4}$ is a skeleton of the same manifold
$(S^2,(2,1),(3,1),(i,1),(1,-1))$ with one vertex less.

\end{document}